\documentclass[12pt]{article}

\usepackage{a4wide}
\usepackage{amssymb}
\usepackage{amsfonts}
\usepackage{amsmath}
\input xy
\xyoption{arrow} \xyoption{matrix}

\date{}

\newtheorem{proposition}{Proposition}[section]
\newtheorem{theorem}[proposition]{Theorem}
\newtheorem{lemma}[proposition]{Lemma}

\newtheorem{corollary}[proposition]{Corollary}

\def\Kdim{{\rm K.dim }\,}

\def\der{\partial }

\def\nFM0{{\nu }_{F,M_0}}
\def\nFN0{{\nu }_{F,N_0}}
\def\nGN0{{\nu }_{G,N_0}}

\def\N0{ {\bf N}_0 }

\def\t{\otimes}
\def\g{\gamma}
\def\v{\varphi}
\def\ra{\rightarrow}

\def\lra{\leftrightarrow}
\def\Xpm{X^{\pm }}

\def\s{\sigma}

\def\l1{{\lambda}_1}

\def\a{\alpha}
\def\a0{ {\alpha }_0}
\def\a1{ {\alpha }_1}

\def\l{\lambda}

\def\nFGM0{{\nu }_{F,G,M_0}}

\def\nFN0{{\nu}_{F,N_0}}

\def\sm{{\sigma}^m}

\def\sm1{{\sigma}^{-1}}

\def\smtp1{{\sigma}^{-t+1}}

\def\S1{S^{-1}}

\def\Xpm1{X^{\pm 1}_1}

\def\sPM1{{\sigma }^{\pm 1}}
\def\sMP1{{\sigma }^{\mp 1 }}

\def\d{\delta}

\def\di{{\rm d.ind}}

\def\CA{{\cal A}}

\def\CD{{\cal D}}

\def\Ytm1{Y^{t-1}}
\def\Yim1{Y^{i-1}}

\def\CF{{\cal F}}

\def\Aut{{\rm Aut}}

\def\bA{\overline{A}}
\def\Der{{\rm Der }}
\def\ad{{\rm ad }}
\def\dim{{\rm dim }}

\def\ker{ {\rm ker } }

\def\D{ \Delta }

\def\SL2Z{ {\rm SL}_2({\bf Z}) }

\def\CR{ {\cal R}}

\def\Gp1{ G^{1 , 1 } }
\def\P11{ P^{-1 , 1 } }
\def\Pp1{ P^{1 , 1 } }

\def\nCLsr{{}^\nu\kern-2pt {\cal L}^{\sigma , \rho  }}
\def\nP{{}^\nu \kern-2pt P}
\def\nL{{}^\nu\kern-2pt L}
\def\nLL{{}^\nu\kern-2pt \Lambda}
\def\nPsr{{}^\nu\kern-2pt P^{\sigma , \rho  }}
\def\nLsr{{}^\nu\kern-2pt L^{\sigma , \rho  }}
\def\nuCL{{}^\nu\kern-2pt  {\cal L}}
\def\nCLsr{{}^\nu\kern-2pt {\cal L}^{\sigma , \rho  }}
\def\nCL1m{{}^\nu\kern-2pt {\cal L}^{-1 , 1  }}
\def\x1nu{x^\frac{1}{\nu}}
\def\xm1nu{x^{-\frac{1}{\nu}}}

\def\tK{\widetilde{K}}

\def\CR{ {\cal R}}

\def\ra{\rightarrow }

\def\nAM0{{\nu }_{{\cal A},M_0}}
\def\nAN0{{\nu }_{{\cal A},N_0}}

\def\Kdim{ {\rm Kdim } }

\def\Der{ {\rm Der }}

\def\CR{ {\cal R }}

\def\det{ {\rm det }}
\def\ad{ {\rm ad }}

\def\bx{\overline{x}}

\def\gm{\mathfrak{m}}

\def\di!{\frac{\der^i}{i!}}
\def\dik!{\frac{\der^k_i}{k!}}
\def\hA{\widehat{A}}
\def\hADi{\widehat{A}^{D_i}}
\def\Ns{\mathbb{N}^s}
\def\grA{{\rm gr}(A)}
\def\grphi{{\rm gr}(\phi)}
\def\grpsi{{\rm gr}(\psi)}
\def\zdai{z_{d , \alpha , i}}
\def\tdai{t_{d , \alpha , i}}
\def\ord{{\rm ord}}

\def\bd{\overline{\delta}}
\def\Adhi{A^{\delta}_{\widehat{i}}}
\def\Adhj{A^{\delta}_{\widehat{j}}}
\def\tK{\widetilde{K}}

\begin{document}

\author{V. V. \  Bavula 
}

\title{Generators and defining relations for ring of invariants of commuting locally nilpotent
derivations or automorphisms}

\maketitle
\begin{abstract}
Let $A$ be an algebra over a field $K$ of characteristic zero, let
$\d_1, \ldots , \d_s\in \Der_K(A)$ be {\em commuting locally
nilpotent}
 $K$-derivations such that $\d_i(x_j)=\d_{ij}$, the Kronecker
 delta, for some elements $x_1,\ldots , x_s\in A$. A set of
 algebra generators for the algebra $A^\d:= \cap_{i=1}^s\ker
 (\d_i)$ is found {\em explicitly} and a set of {\em defining
 relations} for the algebra $A^\d$ is described. Similarly, given
 a set $\s_1, \ldots , \s_s\in \Aut_K(A)$ of {\em commuting}
 $K$-automorphisms of the algebra $A$ such that the maps
 $\s_i-{\rm id_A}$ are {\em locally nilpotent} and $\s_i
 (x_j)=x_j+\d_{ij}$, for some elements $x_1,\ldots , x_s\in A$. A set of algebra generators for the algebra
 $A^\s:=\{ a\in A\, | \, \s_1(a)=\cdots =\s_s(a)=a\}$ is found
 {\em explicitly} and a set of defining relations for the algebra
 $A^\s$ is described. In general, even for a {\em finitely generated noncommutative}
 algebra $A$
  the algebras of invariants $A^\d $ and $A^\s $  are {\em
 not} finitely generated, {\em not} (left or right) Noetherian
  and {\em does not} satisfy finitely many defining
 relations (see examples). Though, for a {\em finitely generated commutative} algebra $A$ {\em
 always}
 the {\em opposite} is true.  The derivations (or automorphisms) just
 described appear often in may different situations after
 (possibly)  a localization of the algebra $A$.

 {\em Mathematics subject classification
2000:  16W22, 13N15, 14R10,  16S15, 16D30.}

$${\bf Contents}$$
\begin{enumerate}
\item Introduction. \item Generators and defining relations for
ring of invariants of commuting locally nilpotent derivations.
\item Generators and defining relations for ring of invariants of
commuting automorphisms. \item The inverse map for  automorphism
that preserve the ring of invariants of derivations. \item
Integral closure and commuting locally nilpotent derivations.
\item A construction of simple algebras. \item Linear maps as
differential operators.
\end{enumerate}
\end{abstract}


\section{Introduction}

The following notation will remain {\bf fixed} throughout the
paper (if it is not stated otherwise): $K$ is a field of
characteristic zero (not necessarily algebraically closed), $A$ is
an (associative, not necessarily commutative) algebra with $1$,
module means a left module, $\d_1, \ldots , \d_s\in \Der_K(A)$ are
{\em commuting locally nilpotent}
 $K$-derivations of the algebra $A$,  $A^\d:= \cap_{i=1}^s\ker
 (\d_i)$ is {\em the algebra of invariants} (or {\em constants}) for the
 derivations $\d := (\d_1, \ldots , \d_s)$; $\s_1, \ldots , \s_s\in \Aut_K(A)$ are {\em
 commuting}
 $K$-automorphisms of the algebra $A$ such that the maps
 $\s_i-{\rm id_A}$ are {\em locally nilpotent} (for each $a\in A$,
 $(\s_i-{\rm id_A})^n(a)=0$ for all $n\gg 1$),  $A^\s:=\{ a\in A\, | \, \s_1(a)=\cdots =\s_s(a)=a\}$
 is {\em the algebra of invariants} for the automorphisms $\s := (\s_1,
 \ldots , \s_n)$.

 Theorem \ref{18Dec05} describes algebras $A$ for which there
 exist a set of commuting locally nilpotent derivations $\d_1,
 \ldots , \d_s$ such that  $\d_i(x_j)=\d_{ij}$, the Kronecker
 delta, for some elements $x_1,\ldots , x_s\in A$ and all $1\leq i,j\leq s$
 (the algebras $A$ are {\em iterated Ore extensions} of a very special
 type). Similarly, Theorem \ref{20Dec05} describes algebras $A$
 for which there exists a set of commuting automorphisms
 $\s_1,\ldots , \s_s$ and  a set of elements $x_1, \ldots ,
 x_s\in A$ such that the maps $\s_i-{\rm id_A}$ are locally
 nilpotent and $\s_i  (x_j)=x_j+\d_{ij}$ for all $1\leq i,j\leq s$
 where ${\rm id}_A$ is the identity map of $A$. The algebras $A$
 are precisely of the type as in Theorem \ref{18Dec05} and vice
 versa. In particular, the problems of finding generators and
 defining relations for the algebra $A^\d$ is the `same' as the
 identical one for the algebra $A^\s$. So, we will restrict
 ourselves mainly to the case of derivations.

{\it Remark}. Two old open problems, the {\em Jacobian Conjecture}
and the {\em Dixmier Problem}, are essentially questions about
 whether certain {\em commuting} derivations $\d_1, \ldots , \d_s$
(of the polynomial algebra or the Weyl algebra, respectively) such
that $\d_i(x_j)=\d_{ij}$ for some elements $x_1, \ldots , x_s$ are
{\em locally nilpotent}. In this paper, we will see that this type
of derivations is more common that one may expect. Typically, this
derivations appear after localization of algebra. In order to
study that kind of derivations it is naturally to look at the
locally nilpotent case first.

 Theorem \ref{15Nov05} gives {\em explicitly} a set of algebra
 generators for the algebra $A^\d$ and describes {\em explicitly}
 the set of defining relations for the generators. More one can say
  in  the important special cases, Corollary \ref{ca18Dec05}
 ($A$ is commutative) and Theorem \ref{19Dec05} (if $[x_i, A^\d
 ]\subseteq A^\d$ for all $i=1,\ldots , s$). Plenty of examples are
 considered. A connection with rings of differential operators is
 described (Corollary \ref{22Dec05}). One can produce an example
 of a finitely generated noncommutative algebra $A$ such that
 the algebras $A^\d $ and $A^\s$
 are {\em not} finitely generated, {\em not} left/right Noetherian, and that generators {\em do not}
 satisfy finitely many defining relations (see Section \ref{2GenCLND}).

Theorem \ref{i19Dec05} gives {\em explicitly} a formula for the
inverse of an automorphism of the algebra $A$ that preserves the
ring of invariants $A^\d$. As an application, we deduce the
inverse formula for an automorphism of the $n$'th Weyl algebra
with polynomial coefficients (Theorem \ref{i8Nov05}).

Theorem \ref{24Dec05} describes algebras $A$ that admit a set of
commuting locally nilpotent derivations with left localizable
kernels. As an application of Theorem \ref{24Dec05} of how to find
explicitly the integral closure $\widetilde{K}$ of the field $K$
in the algebra $A$ is given by Corollary \ref{1c24Dec05}. Theorem
\ref{26Dec05} gives a construction of simple algebras coming from
a set of commuting locally nilpotent derivations.

Let $A=A_n\t P_m$ be the $n$'th Weyl algebra with polynomial
coefficients $P_m$ and $\der_1:=\frac{\der}{\der x_1},\ldots ,
\der_s:=\frac{\der}{\der x_s}\in \Der_K(A)$, $s:= 2n+m$,  be
 the formal partial derivatives of the
algebra $A$, it is a set of commuting locally nilpotent
derivations of the algebra $A$.  Theorem \ref{16Dec05} establishes
a natural isomorphism of the algebra ${\rm End}_K(A)$ and the
algebra $A[[\der_1, \ldots , \der_s]]$, and Theorem \ref{15Dec05}
gives a formula (a sort of a `noncommutative' Taylor formula but
for linear maps rather than for series or polynomials) that
represents any $K$-linear map $a:A\ra A$ as a formal series
$a=\sum_{\alpha \in \Ns} a_\alpha \der^\alpha$, $a_\alpha \in A$.
In particular, for any $\s\in \Aut_K(K[x_1, \ldots , x_m])$, $\s =
\sum_{\alpha \in \mathbb{N}^m} \frac{\prod_{i=1}^m(\s
(x_i)-x_i)^{\alpha_i}}{\alpha !}\der_1^{\alpha_1}\cdots
\der_m^{\alpha_m}$ (Theorem \ref{s15Dec05}).


\section{Generators and defining relations for ring of invariants of commuting locally nilpotent derivations}
\label{2GenCLND}

Let $A$ be an algebra over a field $K$ and let $\d $ be a
 $K$-derivation of the algebra $A$. The kernel $A^\d:=\ker \, \d $ of $\d $ is a
 subalgebra  of $A$, so-called, the {\em algebra of invariants} (or {\em constants}) of $\d $,  the union of the
 vector spaces $N:=N(\d ,A)=\cup_{i\geq 0}\, N_i$ is a positively
 {\em filtered} algebra ($N_iN_j\subseteq N_{i+j}$ for
 all $i,j\geq 0$) where $N_i:= \ker (\d^{i+1})$.
 Clearly, $N_0= A^\d$ and $N:=\{ a\in A \, | \ \d^n (a)=0$
  for some natural $n\}$. A $K$-derivation $\d $ of
 the algebra $A$ is a {\em locally nilpotent } derivation if for
 each element $a\in A$ there exists a natural number $n=n(a)$ such
 that $\d^n(a)=0$. A $K$-derivation $\d $ is locally nilpotent iff
 $A=N(\d , A)$.

Given a ring $R$ and its derivation $d$. The {\em Ore extension}
$R[x;d]$ of $R$ is a ring freely generated over $R$ by $x$ subject
to the defining relations: $xr=rx+d(r)$ for all $r\in R$.
$R[x;d]=\oplus_{i\geq 0}Rx^i=\oplus_{i\geq 0}x^iR$ is a  left and
right free $R$-module. Given $r\in R$, a derivation $(\ad \,
r)(s):=[r,s]=rs-sr$ of $R$ is called an {\em inner} derivation of
$R$.

\begin{lemma}\label{dx=1}
\cite{inform'05} Let $A$ be an algebra over a field $K$ of
characteristic zero and $\d $ be a $K$-derivation of $A$ such that
$\d (x)=1$ for some $x\in A$. Then $N(\d ,A)=A^\d [x; d]$ is the
Ore extension with coefficients from the algebra $A^\d$, and the
derivation $d$ of the algebra $ A^\d$ is the restriction of the
inner derivation $\ad \, x $ of the algebra $A$ to its subalgebra
$A^\d$. For each $n\geq 0$, $N_n=\oplus_{i=0}^n\, A^\d
x^i=\oplus_{i=0}^n\, x^iA^\d$.
\end{lemma}

When the algebra $A$ is commutative the result above is old and
well-known.

\begin{theorem}\label{8Nov05}
\cite{inform'05} Let $A$ be an algebra over a field $K$ of
characteristic zero, $\d $ be a locally nilpotent $K$-derivation
of the algebra $A$ such that $\d (x)=1$ for some $x\in A$.  Then
the $K$-linear map $\phi :=\sum_{i\geq 0} (-1)^i\frac{x^i}{i!}\d^i
:A\ra A$ (resp. $\psi :=\sum_{i\geq 0} (-1)^i\d^i (\cdot
)\frac{x^i}{i!} :A\ra A$) satisfies the following properties:
\begin{enumerate}
\item  $\phi$ (resp. $\psi $) is a homomorphism of right (resp.
left) $A^\d$-modules. \item $\phi$ (resp. $\psi $) is a projection
onto the algebra $A^\d$:
\begin{eqnarray*}
\phi : & A=A^\d \oplus xA\ra A^\d \oplus xA, \;\; a+xb\mapsto a,
\;\; {\rm where}\;\; a\in A^\d, \; b\in A,  \\
 \psi : & A=A^\d \oplus Ax\ra A^\d \oplus Ax, \;\; a+bx\mapsto a,
\;\; {\rm where}\;\; a\in A^\d, \; b\in A.
\end{eqnarray*}
 In particular,  ${\rm
im} (\phi )={\rm im} (\psi )= A^\d$ and $\phi (y)=y=\psi (y)$ for
all $y\in A^\d$. \item $\phi (x^i)=\psi (x^i)=0$, $i\geq 1$. \item
$\phi$ and $\psi$ are algebra homomorphisms provided $x\in Z(A)$,
the centre of the algebra $A$.
\end{enumerate}
\end{theorem}

The following notation will remain {\em fixed} till the end of
this section (if it is not stated otherwise): $A$ is an algebra
over a field $K$ of characteristic zero,
 $\d_1, \ldots , \d_s$ are {\em commuting locally nilpotent}
 $K$-derivations of $A$,  $A^\d
 := \cap_{i=1}^s A^{\d_i}$ is the {\em algebra  of invariants} for
 the set of derivation $\d_1, \ldots , \d_s$
 where $A^{\d_i}:= \ker (\d_i)$. The algebra $A$ is equipped with
 the filtration $\{ N_i\}_{i\geq 0}$  ($N_iN_j\subseteq
 N_{i+j}$, for  $i,j\geq 0$) where $N_i:= \{ a\in A\, | \,
 \d^\alpha (a)=0$ for all $\alpha =(\alpha_i) \in \Ns $ with
 $|\alpha |:= \alpha_1+\cdots +\alpha_s>i\}$,  where $\d^\alpha := \d_1^{\alpha_1}\cdots
 \d_s^{\alpha_s}$. $A=\cup_{i\geq 0}N_i$, $N_0:= A^\d \subset
 N_1\subset \cdots$. For $0\neq a\in A$, a unique number  $i$ such that $a\in
 N_i\backslash N_{i-1}$ is called the {\em order} of $a$, denoted $\ord
 (a)$. Consider the associated graded algebra ${\rm gr} (A):= \oplus_{i\geq
 0}N_i/N_{i-1}$ ($N_{-1}:=0$).

The next theorem is a crucial step  in many results that follow.
\begin{theorem}\label{18Dec05}
Let $A$ be an arbitrary algebra over a field $K$ of characteristic
zero. The following statements are equivalent.
\begin{enumerate}
\item There exist commuting locally nilpotent $K$-derivations
$\d_1, \ldots , \d_s$ of the algebra $A$ and elements $x_1,
\ldots, x_s\in A$ satisfying $\d_i(x_j)=\d_{ij}$, the Kronecker
delta. \item The algebra $A$ is an iterated Ore extension $A=
B[x_1; d_1]\cdots [ x_s; d_s]$ such that $d_i (B)\subseteq B$ and
$d_i(x_j)\in B$ for all $1\leq i,j\leq s$.
\end{enumerate}
If, say, the first condition holds, then $A=A^\d [x_1; d_1] \cdots
[ x_s; d_s]$ is an iterated Ore extension of the ring of
invariants $A^\d :=\cap_{i=1}^s A^{\d_i}$ such that $d_i := \ad
(x_i)$, $[x_i, A^\d ] \subseteq A^\d$, and $[x_i, x_j]\in A^\d$
for all $i,j$. In particular, $A=\oplus_{\alpha \in \Ns }x^\alpha
A^\d = \oplus_{\alpha \in \Ns } A^\d x^\alpha$ where $x^\alpha :=
x_1^{\alpha_1}\cdots x_s^{\alpha_s}$, and $A=\cup_{i\geq 0} N_i$
where $N_i=\oplus_{|\alpha |\leq i}x^\alpha A^\d = \oplus_{|\alpha
|\leq i} A^\d x^\alpha$ for $i\geq 0$.
\end{theorem}

{\it Proof}. $(1\Rightarrow 2)$ Applying Lemma \ref{dx=1} step by
step we have the result ($A=\oplus_{k\geq 0}A^{\d_s}x_s^k$, if
$i<s$ then $x_i=\sum \l_{ij}x_s^j$ for some $\l_{ij}\in A^{\d_s}$;
now $0=\d_s (x_i)=\sum j\l_{ij}x_s^{j-1}$ implies $x_i\in
A^{\d_s}$): 
\begin{equation}\label{AAds}
A=A^{\d_s}[x_s; d_s]=(A^{\d_{s-1}}\cap A^{\d_s}) [x_{s-1};
d_{s-1}] [x_s; d_s]=\cdots =A^\d [x_1; d_1]\cdots [ x_s; d_s],
\end{equation}
where  $d_i:= \ad (x_i)$. For all $i,j,k$,  $ \d_k ([x_i,
x_j])=\d_{ki}[1,x_j]+\d_{kj}[x_i, 1]=0$ and $\d_k([x_i, A^\d
])=\d_{ki}[1,A^\d ]=0$, hence all $[x_i, x_j]\in A^\d$ and $[x_i,
 A^\d ]\subseteq A^\d$.

$(2\Rightarrow 1)$  Given an algebra $A$ as in the second
statement. The formal partial derivatives $\frac{\der}{\der x_1},
\ldots ,\frac{\der}{\der x_s}\in \Der_B(A)$ satisfy the condition
of the first statement.  $\Box $

It is obvious that the elements $x_1, \ldots , x_s$ are {\em not}
(left and right) zero divisors in $A$. Next, we have many examples
of derivations as in Theorem \ref{18Dec05}.

{\it Example}. Let $F_n:=K\langle x_1, \ldots , x_n\rangle$ be a
free algebra over the field $K$, $\der_1:=\frac{\der}{\der x_1},
\ldots ,\der_n:=\frac{\der}{\der x_n}\in \Der_K(F_n)$ be the
formal  partial derivatives and $I$ be an ideal of $F_n$ which is
$\der$-invariant (that is $\der_i(I)\subseteq I$ for all $i$). The
induced derivations $\d_1, \ldots \d_n\in \Der_K(A)$ where
$A:=F_n/I$, $\d_i (f+I)=\der_i(f)+I$, $f\in F_n$, are commuting
locally nilpotent derivations of the algebra $A$ and $\d_i
(\overline{x}_i)=\d_{ij}$ for all $1\leq i,j\leq n$ where
$\bx_i:=x_i+I$. If the ideal $I$ is generated by the commutators
$[x_i, x_j]$, $1,\leq i,j\leq n$, we have a polynomial algebra
$K[\bx_1, \ldots , \bx_n ]$ and the derivations
$\d_1:=\frac{\der}{\der \bx_1}, \ldots ,\d_n:=\frac{\der}{\der
\bx_n}$.

\begin{corollary}\label{f26Dec05}
Let $A$, $\d_1, \ldots , \d_s$ and $x_1, \ldots , x_s$ be as in
Theorem \ref{18Dec05}, and $\gm $ be a (two sided) ideal of the
algebra $A^\d$ which is $\ad (x_i)$-invariant for all $i=1,\ldots
, s$, and $(\gm ):= A\gm A$ be the ideal of the algebra $A$
 generated by $\gm$, and $A\ra \bA := A/(\gm )$, $a\mapsto
 \overline{a}:= a+(\gm )$. Then $\bA , \bd_1, \ldots , \bd_s$ and
 $\bx_1, \ldots , \bx_s$ satisfy the conditions of Theorem
 \ref{18Dec05}  (where $\bd_i\in \Der_K(\bA )$, $\overline{a}\mapsto
 \overline{\d_i(a)}$), $\bA^{\bd }=\overline{A^\d }$, and
 $N'_i=\overline{N}_i=\oplus_{|\alpha |\leq i}\overline{A^\d }\bx^\alpha =\oplus_{|\alpha |\leq i}\bx^\alpha
 \overline{A^\d }$, for $i\geq 0$, where $\{ N_i\}$ and $\{ N_i'\}$ are the
 filtrations of the algebras $\bA $ and $A$ respectively.
\end{corollary}

{\it Proof}. The derivations $\bd_1, \ldots , \bd_s$ of the
algebra $\bA$ are commuting locally nilpotent derivations such
that $\bd_i(\bx_j)=\d_{ij}\bx_i$, hence they satisfy the
conditions of Theorem \ref{18Dec05}. In particular, $\bA
=\bA^{\overline{\d}} [\bx_1; \overline{d}_1]\cdots [\bx_s;
\overline{d}_s]=\oplus_{\alpha \in \Ns}\bA^{\overline{\d}}
\bx^\alpha =\oplus_{\alpha \in \Ns} \bx^\alpha
\bA^{\overline{\d}}$ and $N_i'=\oplus_{|\alpha |\leq
i}\bA^{\overline{\d}} \bx^\alpha$, $i\geq 0$.
 On the other hand, $A=\oplus_{\alpha \in \Ns}A^\d x^\alpha$ and
 $(\gm )=\oplus_{\alpha \in \Ns}\gm x^\alpha$, hence $\bA =A/(\gm
 )=\oplus_{\alpha \in \Ns}(A^\d /\gm  )\bx^\alpha$. Comparing the
 two direct sums for $\bA$ we must have $\bA^{\overline{\d}} =\overline{A^\d }$, and
 $N'_i=\overline{N}_i=\oplus_{|\alpha |\leq i}\overline{A^\d
 }\bx^\alpha =\oplus_{|\alpha |\leq i}\bx^\alpha \overline{A^\d
 }$, for $i\geq 0$.
$\Box $

The next result is a criterion of when the ring of invariants
$A^\d$ is left/right Noetherian.

\begin{corollary}\label{c23Dec05}
Let $A$, $\d_1, \ldots , \d_s$ and $x_1, \ldots , x_s$ be as in
Theorem \ref{18Dec05}.  Then the following statements are
equivalent:
\begin{enumerate}
\item The algebra $A$ is left (resp. right) Noetherian. \item The
algebra $A^\d$ is  left (resp. right) Noetherian.\item The algebra
$\grA $ is  left (resp. right) Noetherian.
\end{enumerate}
\end{corollary}

{\it Proof}. $(1\Leftrightarrow 2)$ It is a well-known fact that
if a coefficient ring is left (resp. right) Noetherian then so is
an  iterated Ore extension, and vice versa (use iteratively an
analogue of the Hilbert Basis Theorem for Ore extensions). Now,
the first two statements are equivalent by Theorem \ref{18Dec05}.

$(2\Leftrightarrow 3)$ The associated graded algebra ${\rm gr}
(A)\simeq A^\d [ \bx_1; \overline{d}_1]\cdots  [ \bx_s;
\overline{d}_s]$ is an iterated Ore extension where $\bx_i:=
x_i+N_1/N_0$, $\overline{d}_i=\ad (\bx_i)$, and $[\bx_i, \bx_j]=0$
for all $i,j$. Now, repeat the above argument. $\Box $

{\it Remark}. Using the previous proof one can write down several
similar statements for properties that are `stable' under the
operations of taking iterated Ore extension and ${\rm gr} (\cdot
)$ ({\it eg}, `being domain', etc). For a property of `being
finitely generated algebra', in general, it is not true that `$A$
is finitely generated  $\Rightarrow $ $A^\d $ is  finitely
generated' (see an example after Theorem \ref{15Nov05}), but for
commutative algebras it is the case (Corollary \ref{1c18Dec05}).

Corollary \ref{22Dec05} provides natural examples of commuting
locally nilpotent derivations (on non-commutative algebras), it
also shows that the order filtration $\{ \CD (R)_i\}$ on the ring
of differential operators is, in fact, the filtration $\{ N_i\}$
for certain commuting locally nilpotent derivations of $\CD (R)$
(this fact may simplify arguments in finding {\em explicitly} the
ring $\CD (R)$ of differential operators in certain cases, see the
example below).

Let $R$ be a commutative finitely generated $K$-algebra and $\CD
(R)=\cup_{i\geq 0}\CD (R)_i$ be its {\em ring of differential
operators} on the ring $R$ equipped with the {\em order
filtration} $\{ \CD (R)_i\}$: $\CD (R)$ is a $K$-subalgebra of the
algebra ${\rm End}_K(A)$ where $\CD (R)_0:= {\rm End}_R(R)\simeq
R$ $((x\mapsto rx)\lra r)$, and
$$ \CD (R)_i:= \{ u \in {\rm End}_K(A): \;\; [u,r]\in \CD
(R)_{i-1}\;\; {\rm for \;\; all}\;\; r\in R\}.$$

\begin{corollary}\label{22Dec05}
Let a domain $R$ be a commutative finitely generated $K$-algebra
of Krull dimension $n>0$, $\CD (R)=\cup_{i\geq 0}\CD (R)_i$ be the
ring of differential operators on $R$, $x_1, \ldots , x_n$ be
algebraically independent (over $K$)  elements of $R$. Then
\begin{enumerate}
\item $\d_1:= \ad (x_1), \ldots , \d_n:=\ad (x_n)$ is a set of
commuting locally nilpotent derivations of the algebra $\CD (R)$.
\item The order filtration $\{ \CD (R)_i\}$ coincides with the
filtration $\{ N_i\}$ associated with the derivations $\d_1,
\ldots , \d_n$, i.e. $\CD (R)_i=N_i$ for all $i\geq 0$. In
particular, $\CD (R)^\d =R$.
\end{enumerate}
\end{corollary}

{\it Proof}. The first statement is obvious. To prove the second
statement, note that $\CD (R)_i\subseteq N_i$ for all $i\geq 0$
which follows directly from the definitions of both filtrations.

Let $P_n:= K[x_1, \ldots , x_n]$ and $Q_n=K(x_1, \ldots , x_n)$ be
its field of fractions. The field $Q:={\rm Frac} (R)$ of fractions
of $R$ is a {\em finite separable} field extension of $Q_n$. It
well-known that one can pick up a nonzero element, say $r\in R$,
such that the localization $R_r:= R[r^{-1}]$ of $R$ at the powers
of the element $r$ is a {\em regular} domain, $\der_i
(R_r)\subseteq R_r$ for all $i$, and
$\Der_K(R_r)=\oplus_{i=1}^nR_r\der_i$  where $\der_i:=
\frac{\der}{\der x_i}$ are the partial derivatives of $Q_n$
uniquely extended to derivations of the field $Q$. Since the
algebra $R_r$ is regular the ring of differential operators $\CD
(R_r)$ on the algebra $R_r$ is generated by the algebra $R_r$ and
$\Der_K(R_r)$, hence $\CD (R_r)=\oplus_{\alpha \in \mathbb{N}^n}
R_r\der^\alpha$ and $\CD (R_r)_i=\oplus_{|\alpha |\leq i}
R_r\der^\alpha$, $i\geq 0$. Comparing these equalities with
similar ones from Theorem \ref{18Dec05}: $\CD (R_r)=\oplus_{\alpha
\in \mathbb{N}^n} \CD (R_r)^\d \der^\alpha =\cup_{i\geq
0}N_i(R_r)$ and $N_i(R_r)=\oplus_{|\alpha |\leq i } \CD (R_r)^\d
\der^\alpha$, $i\geq 0$, and taking into account the inclusions
$\CD (R_r)_i\subseteq N_i(R_r)$ for $i\geq 0$,  we must have $\CD
(R_r)^\d = R_r$ and $N_i (R_r)= \CD (R_r)_i$, $i\geq 0$. Since
$\CD (R)\subseteq \CD (R_r)$ and $\CD (R)_i= \CD (R)\cap \CD
(R_r)_i$ for all $i\geq 0$, and $N_i=\CD (R)\cap N_i(R_r)$, $i\geq
0$, we conclude that $N_i=\CD (R)_i$ for all $i\geq 0$.  $\Box $

{\it Example}. As an application of Theorem \ref{18Dec05} and
Corollary \ref{22Dec05}, let us give a short proof of the
well-known fact that {\em the ring of differential operators $\CD
(P_n)$ on a polynomial algebra $P_n:=K[x_1, \ldots , x_n]$
(so-called, the Weyl algebra) is generated by $P_n$ and the
partial derivatives} $\der_1, \ldots , \der_n$ of $P_n$: the inner
derivations $\d_1:= -\ad (x_1) , \ldots , \d_n:=-\ad (x_n)$ of the
algebra $E:={\rm End}_K(P_n)$ commute. Let $N$ be the {\em
largest} subalgebra of $E$ on which all the derivations $\d_i$ act
locally nilpotently (take the sum of all the subalgebras of $E$
with the last property).  Clearly, $\CD (P_n)\subseteq N$ and
$N^\d =E^\d = {\rm End}_{P_n}(P_n)\simeq P_n$. By Theorem
\ref{18Dec05}, $N=P_n\langle \der_1, \ldots , \der_n\rangle
\subseteq \CD (P_n)$, hence $N=\CD (P_n)$, and, by Corollary
\ref{22Dec05}, $\CD (P_n)_i=N_i=\oplus_{|\alpha | \leq
i}P_n\der^\alpha$ for all $i\geq 0$.

\begin{corollary}\label{1c18Dec05}
Let $A$, $\d_1, \ldots , \d_s$ and $x_1, \ldots , x_s$ be as in
Theorem \ref{18Dec05}. Suppose that the elements $x_1, \ldots ,
x_s$ are central. Then the following statements are equivalent:
\begin{enumerate}
\item The algebra $A$ is finitely generated. \item The algebra
$A^\d$ is finitely generated.\item The algebra $\grA $ is finitely
generated.
\end{enumerate}
\end{corollary}

{\it Proof}.  Since the elements $x_1, \ldots , x_s$ are central,
by Theorem \ref{18Dec05}, $A\simeq A^\d [x_1, \ldots , x_s]\simeq
\grA$ and $A^\d \simeq A/(x_1, \ldots , x_s)$. Now, it is obvious
that the statements are equivalent.  $\Box $

{\em Till the end of this section} we will {\em assume} that for
the commuting locally nilpotent derivations $\d_1, \ldots , \d_s$
of $A$ there exist elements  $x_1, \ldots , x_s\in A$ such that
$\d_i(x_j)=\d_{ij}$, the  Kronecker delta.

For each $i=1, \ldots , s$, consider the maps from Theorem
\ref{8Nov05},
$$ \phi_i:=\sum_{k\geq 0}(-1)^k\frac{x_i^k}{k!}\d_i^k,  \;\; \psi_i:=\sum_{k\geq 0}(-1)^k\d_i^k
(\cdot ) \frac{x_i^k}{k!}\, : \, A\ra A.$$ The maps $\phi_i$ and
$\psi_i$ are homomorphisms of {\em right} and {\em left}
$A^\d$-modules respectively. The maps 
\begin{equation}\label{rphis}
\phi:=\phi_s\phi_{s-1}\cdots \phi_1:A\ra A, \;\; a=\sum_{\alpha
\in \mathbb{N}^s}  x^\alpha\l_\alpha\mapsto \phi (a)=\l_0,
\end{equation}
\begin{equation}\label{lphis}
\psi:=\psi_1\psi_2\cdots \psi_s:A\ra A, \;\; a=\sum_{\alpha \in
\mathbb{N}^s} \l_\alpha x^\alpha\mapsto \psi (a)=\l_0,
\end{equation}
are {\em projections} onto the subalgebra $A^\d$ of $A=A^\d \oplus
(\oplus_{0\neq \alpha \in \Ns }x^\alpha A^\d )$ and $A= A^\d
\oplus (\oplus_{0\neq \alpha \in \Ns } A^\d x^\alpha )$
respectively, they are homomorphisms of right and left $A^\d
$-modules respectively.

\begin{theorem}\label{a18Dec05}
Let $A$ be as in Theorem \ref{18Dec05}.  For any $a\in A$,
$$ a=\sum_{\alpha \in \mathbb{N}^s}x^\alpha \phi (\frac{\d^\alpha}{\alpha
!} a)=\sum_{\alpha \in \mathbb{N}^s}\psi (\frac{\d^\alpha}{\alpha
!} a)x^\alpha.$$
\end{theorem}

{\it Proof}. If $a=\sum  x^\alpha \l_\alpha $, $\l_\alpha \in
A^\d$, then, by (\ref{rphis}), $\phi (\frac{\d^\alpha}{\alpha !}
a)=\l_\alpha$. Similarly, if $a=\sum  \l_\alpha x^\alpha  $,
$\l_\alpha \in A^\d$, then, by (\ref{lphis}), $\psi
(\frac{\d^\alpha}{\alpha !} a)=\l_\alpha$. $\Box $

So, the identity map ${\rm id} : A \ra A$ has  nice presentations
\begin{equation}\label{ida}
{\rm id}(\cdot ) = \sum_{\alpha \in \mathbb{N}^s}x^\alpha \phi
(\frac{\d^\alpha}{\alpha !} (\cdot ))= \sum_{\alpha \in
\mathbb{N}^s}\psi (\frac{\d^\alpha}{\alpha !} (\cdot ))x^\alpha .
\end{equation}

  Clearly, $\phi (N_i)\subseteq
 N_i$ and $\psi (N_i)\subseteq
 N_i$ for all  $i\geq 0$.  Consider
 the associated graded algebra ${\rm gr} (A):= \oplus_{i\geq
 0}N_i/N_{i-1}$ ($N_{-1}:=0$). So, let $\grphi , \grpsi  : \grA \ra \grA$ be the
 induced maps (for $\overline{a} = a+N_{i-1}\in N_i/ N_{i-1}$, ${\rm gr} (\phi ) (\overline{a} )=\phi (a)+N_{i-1}$
 and ${\rm gr} (\psi ) (\overline{a} )=\psi (a)+N_{i-1}$). Let $D$ be
  a free
 multiplicative monoid
 generated freely by the inner derivations $\ad (x_1), \ldots , \ad (x_s)$ of the algebra $A$.
  There is an obvious action of $D$ on the algebra $A$ (and an obvious linear map $D\ra {\rm End}_K(A)$). Let $\{ y_i\, | \, i\in
 I\}$ be a set of algebra generators for $A$.  For
 each $d\in D$ and $i\in I$, let $\zdai := d\phi (\frac{\d^\alpha }{\alpha!}
 y_i)$ where $\alpha !:= \alpha_1! \cdots \alpha_s!$ and $\alpha_i
 \leq \ord (y_i)$ for all $i=1, \ldots , s$. Let ${\rm id}_A$ be the identity map of $A$.   For
 each $d'\in D^*:= D\backslash \{ {\rm id}_A \}$ and $j=1, \ldots , s$, let $x_{d',j}:= d' (x_j)$.
  By Theorem \ref{18Dec05}, all the elements $\zdai , x_{d',j}\in A^\d$.  For each $\zdai$ and each $x_{d',j}$ we attach
 (noncommutative)
 variables
 $\tdai$ and $X_{d',j}$ respectively. Let $\CF := K\langle \tdai , X_{d',j} \, |\, d\in D, i\in I, \alpha , d'\in D^*,  1\leq j \leq s\rangle$ be
 a free associative algebra and $f(\tdai , X_{d',j})=f( \{ \tdai , X_{d',j} \, |\, d\in D, i\in I, \alpha ,  d'\in D^*, 1\leq j \leq s\}) $
 be a typical element
 of $\CF$ (the symbols in the brackets, i.e. $\tdai , X_{d',j}$, stand for {\em all}
 the non-commutative arguments of the element $f$).

\begin{theorem}\label{15Nov05}
The algebra $A^\d$ is generated by all the elements $\{ \zdai ,
x_{d', j}\}$ that satisfy the defining relations $\CR =\{ f(\tdai
, X_{d', j} )\in \CF \, | \, f(\zdai , x_{d', j}  )\in
\sum_{i=1}^s x_iA\}$. Similarly, the algebra $A^\d$ is generated
by all the elements $\{ \zdai':= d\psi
(\frac{\d^\alpha}{\alpha!}y_i), x_{d', j}\}$ that satisfy the
defining relations $\CR' =\{ f(\tdai , X_{d', j} )\in \CF \, | \,
f(\zdai' , x_{d', j}  )\in \sum_{i=1}^s Ax_i\}$.
\end{theorem}

{\it Proof}. Recall that $A=\oplus_{\alpha \in \Ns } x^\alpha
A^\d$, and so each $y_i$ is a unique sum $y_i= \sum x^\alpha
y_{\alpha , i}$ where $y_{\alpha , i}=\phi
(\frac{\d^\alpha}{\alpha !}y_i)\in A^\d$ (Theorem \ref{a18Dec05})
and $\alpha_i \leq \ord (y_i)$ for all $i=1, \ldots , s$
(otherwise, $y_{\alpha , i}=0$). The set $\{ y_i\, | \, i\in I\}$
is a set of $K$-algebra generators for $A$, hence so is the set
$\{ y_{\alpha , i}, x_1, \ldots , x_s\, | \, i\in I, \alpha \}$
(with obvious restrictions on $\alpha \in \Ns$ for each $y_{\alpha
, i}$, that is $\alpha_i \leq \ord (y_i)$ for all $i=1, \ldots ,
s$). Since all the $y_{\alpha , i}\in A^\d$, $[x_j, x_k]\in A^\d$,
and $[x_j, A^\d ]\subseteq A^\d $ for all $j,k$, any element $a\in
A$ can be written as a sum $a=\sum x^\alpha a_\alpha$ where each
coefficient $a_\alpha $ belongs to the subalgebra, say $\CA$, of
$A$ generated by all the elements $\{ \zdai, x_{d', j}\}$ in the
theorem. It follows that $A^\d =\phi (A)\subseteq \CA$, the
opposite inclusion, $\CA \subseteq A^\d$, is obvious. Therefore,
$A^\d =\CA$.

Since all the elements $\zdai , x_{d', j}\in A^\d$ and the map
$\phi $ is a projection onto the ring of invariants $A^\d$, an
element $f(\tdai , X_{d', j})\in \CF$ is a relation for the set of
generators $\{ \zdai, X_{d', j}\}$ of the algebra $A^\d$, i.e.
$f(\zdai, x_{d', j})=0$, iff $\phi (f(\zdai, x_{d', j}))=0$ iff
$f(\zdai, x_{d', j})\in \sum_{j=1}^sx_jA$.

To prove the remaining case, repeat the above arguments making
obvious adjustments.  $\Box $

\begin{corollary}\label{c15Nov05}
Let $\{ a_j\, | \, j\in J\}$ be algebra generators for $A^\d$ and
$F_J=K\langle Y_j\, | \,  j\in J\rangle $ be a free algebra. Then
$\CR =\{ f(Y_j)\in F_J\, | \, f(a_j)\in \sum_{i=1}^s x_iA\}$
(resp. $\CR' =\{ f(Y_j)\in F_J\, | \, f(a_j)\in \sum_{i=1}^s
Ax_i\}$) are defining relations for the algebra $A^\d$.
\end{corollary}

{\it Proof}. Repeat the arguments as in the proof of Theorem
\ref{15Nov05}.  $\Box $

{\it Remark}. The proof of Theorem \ref{15Nov05} shows that the
choice of generators there might be not the most economical one if
the algebra $A$ is far from being free (see also Corollary
\ref{ca18Dec05} and Theorem \ref{19Dec05}). The proof of Theorem
\ref{15Nov05} shows that in order to find generators and defining
relations for the algebra $A^\d$ one should
\begin{enumerate}
\item take algebra generators $\{ y_i\}_{i\in I}$ for the algebra
$A$, \item find the coefficients $y_{\alpha , i}:= \phi
(\frac{\d^\alpha}{\alpha !}y_i)\in A^\d$ of each element $y_i=\sum
x^\alpha y_{\alpha , i}$, \item choose a basis, say $\{ b_j\, | \,
j\in I'\}$, of the $K$-linear span of all the coefficients $\{
y_{\alpha , i}\}$, \item then the algebra $A^\d$ is generated by
the elements $\{ d (b_j), x_{d', i}\, | \, d\in D, d'\in D^*, j\in
I', i=1, \ldots , s\}$, \item choose more economically (if
possible) an algebra generators say $\{ a_j\, | \, j\in J\}$ for
$A^\d$, \item Corollary \ref{c15Nov05} gives the defining
relations.
\end{enumerate}

{\it Example}. Let $A=F_n=\langle x_1, \ldots , x_n\rangle$ be a
free associative algebra over $K$, $\d_1:= \frac{\der}{\der x_1},
\ldots ,\d_s:= \frac{\der}{\der x_s}\in \Der_K(F_n)$ be formal
partial derivatives, $s\leq n$: $\d_i(x_j)=\d_{ij}$, $1\leq i\leq
s$, $1\leq j\leq n$. Then $\phi (x_1)=\cdots = \phi (x_s)=0$ and
$\phi (x_{s+1})=x_{s+1}, \ldots , \phi (x_n)=x_n$. By Theorem
\ref{15Nov05}, $A^\d = K\langle d (x_i), d'(x_j)\, | \,  d\in D,
d'\in D^*,  i=s+1, \ldots , n; j=1, \ldots , s\rangle = K\langle
x_i, d'(x_j)\, | \, d'\in D^*,  i=s+1, \ldots , n; j=1, \ldots ,
s\rangle $, and, by Corollary \ref{c23Dec05}, for $n\geq 2$, the
algebra $A^\d$ is {\em not} left/right Noetherian since $F_n$ is
not. In the special case when $s=n$, the algebra $A^\d$ is a {\em
free algebra in infinitely many variables}. More precisely, it is
generated freely by the elements $\{ (\ad \, x_1)^{\alpha_1}\cdots
(\ad \, x_n)^{\alpha_n}([x_i,x_j])\, | \, i<j, (\alpha_1, \ldots ,
\alpha_n)\in \mathbb{N}^n\}$ (Prop. 2, \cite{Ger98Arch}). Now,
taking any ideal $\gm $ of the algebra $A^\d$ which is {\em not}
finitely generated as an $A^\d$-bimodule and such that the algebra
$A^\d/ \gm $ is not left/right Noetherian, and using Corollary
\ref{f26Dec05} one produces an example of an algebra
$\overline{A}^{\overline{\d}}$ which is {\em not} finitely
generated, {\em not} left/right Noetherian, and {\em does not}
satisfy finitely many defining relations.

\begin{corollary}\label{ca18Dec05}
If, in addition, $A$ is commutative then the map $\phi : A \ra
A^\d$ is an algebra epimorphism with $\ker (\phi )=(x_1, \ldots ,
x_s)$. In particular, $A^\d \simeq A /(x_1, \ldots , x_s)$.
Alternatively, the algebra $A^\d$ is generated by the elements
$\{\phi (y_i)\, | \, i\in I\}$ that satisfy  the defining
relations $\CR = \{ f(t_i )\in K[t_i\, | \, i\in I]  \, | \,
f(\phi (y_i) )\in (x_1, \ldots , x_s)\}$.
\end{corollary}

{\it Proof}. Each of the maps $\phi_i=\psi_i$ is an  algebra
homomorphism, hence so is their product $\phi =\psi$. Now, the
result follows from (\ref{rphis}) or (\ref{lphis}).  $\Box $

{\it Example}. {\em The Weitzenb\"{o}ck derivation} $\d =
x_1\frac{\der }{\der x_2}+x_2\frac{\der }{\der x_3}+\cdots
+x_{n-1}\frac{\der }{\der x_n}$ of the polynomial algebra $P_n:=
K[x_1, \ldots , x_n]$ ($n\geq 3$) is locally nilpotent with $\d
(x_1)=0$. The derivation can be (uniquely) extended to a locally
nilpotent derivation of the localization $P_n[x_1^{-1}]=K[x_1,
x_1^{-1}][ x_2, \ldots , x_n]$ with $\d (x)=1$ where $x:=
\frac{x_2}{x_1}$. By Corollary \ref{ca18Dec05}, the algebra of
invariants $P_n[x_1^{-1}]^\d$ is equal to  $K[x_1, x_1^{-1}][ \phi
(x_3), \ldots , \phi (x_n)]$, the  polynomial algebra in $\phi
(x_i)$, $i=3, \ldots , n$,
  with
coefficients from $K[x_1, x_1^{-1}]$ (note that $\phi
(x_2)=x_1\phi (x)=x_10=0$) where $$\phi (x_i)=\sum_{k=0}^{i-1}
(-1)^k(\frac{x_2}{x_1})^k\frac{x_{i-k}}{k!}, \;\;\ i=3, \ldots ,
n.$$ Hence, $P_n^\d = P_n\cap P_n[x_1^{-1}]^\d =P_n\cap K[x_1,
x_1^{-1}][ \phi (x_3), \ldots , \phi (x_n)]$. Since $\d
(\sum_{i=1}^s Kx_i)\subseteq \sum_{i=1}^s Kx_i$, by the Theorem of
Weitzenb\"{o}ck, $P_n^\d$ is a {\em finitely generated} algebra.
It is an open problem to find explicitly a set of algebra
generators for it (it would imply an explicit description of all
$SL_2$-invariants which is another open problem, in fact, these
two problems are equivalent).

\begin{theorem}\label{19Dec05}
If $[x_i, A^\d ]=0$, $1\leq i \leq s$, then the map $\grphi : \grA
\ra A^\d$ is an algebra epimorphism with kernel generated by the
elements $\bx_1, \ldots , \bx_s$ (where $\bx_j:= x_j+A^\d \in
N_1/N_0$). In particular, $A^\d \simeq \grA /(\bx_1, \ldots ,
\bx_s)$. Alternatively, the algebra $A^\d$ is generated by the
elements $\{\grphi (y_i)\, | \, i\in I\}$ that satisfy  the
defining relations $\CR = \{ f(t_i )\in K\langle t_i\, | \, i\in
I\rangle  \, | \, f(\grphi (y_i) )\in \sum_{j=1}^s \bx_j\grA \}$.
\end{theorem}

{\it Proof}. Since  $[x_i, A^\d ]=0$ and $[x_i, x_j]\in A^\d$ for
all $1\leq i, j \leq s$, it follows from Theorem \ref{18Dec05}
that $\grA \simeq A^\d [ \bx_1, \ldots , \bx_s]$ is a polynomial
algebra over $A^\d$ in $\bx_i:= x_i+A^\d$, $i=1, \ldots , s$. The
induced derivations $\bd_1, \ldots , \bd_s\in \Der_K (\grA )$ of
graded degree $-1$  are commuting locally nilpotent derivations of
the algebra $\grA$ (where $\bd_i: N_j/N_{j-1}\ra N_{j-1}/N_{j-2}$,
$a+N_{j-1}\mapsto
 \d_i (a)+N_{j-2}$) with $\bd_i (\bx_j)=\d_{ij}$. Now, we are in
 the
 situation of Corollary \ref{ca18Dec05}. Let $\overline{\phi}$ be
 the corresponding map from Corollary \ref{ca18Dec05}. Clearly,
  $\overline{\phi} = \grphi $. Now, the result becomes obvious due
  to Corollary \ref{ca18Dec05}. $\Box$

\begin{lemma}\label{r22Dec05}
Let $A$, $\d_1, \ldots , \d_s$, and $x_1, \ldots , x_s$ be as in
Theorem \ref{18Dec05}. If $A'$ is a $\d$-invariant subalgebra of
the algebra $A$ ($\d_i(A')\subseteq A'$ for all $i$) then the
restrictions $\d_1':= \d_1|_{A'}, \ldots ,  \d_s':= \d_s|_{A'}$
are commuting locally nilpotent derivations of the algebra $A'$
and $N_i'=A'\cap N_i$  for all $i\geq 0$, in particular,
$(A')^{\d'}=A'\cap A^\d$, and ${\rm gr}(A')\subseteq {\rm gr}(A)$
is a natural inclusion of graded algebras.
\end{lemma}

{\it Proof}. Obvious. $\Box $

{\it Example}. Given a $K$-algebra $A$ and commuting locally
nilpotent derivations $\d_1, \ldots , \d_s$ of the algebra $A$,
and let $\{ N_i\}$ be the corresponding filtration. Let $Z(A)$ and
$NZD(A)$ be the centre and the set of all the (left and right)
non-zero-divisors of $A$ respectively. Consider the set $S:= A^\d
\cap Z(A)\cap NZD (A)$. The algebra $A$ is a subalgebra of the
localization $\S1 A$ of the algebra $A$ at $S$, the derivations
$\d_1, \ldots , \d_s$ can be uniquely extended to derivations of
the algebra $\S1 A$, denoted in the same fashion. These extended
derivation $\d_1, \ldots , \d_s$ are commuting locally nilpotent
derivations of the algebra $\S1 A$. Suppose that there are
elements $x_1, \ldots , x_s\in \S1 A$ such that $\d_i
(x_j)=\d_{ij}$ for all $i,j$. By Lemma \ref{r22Dec05}, $N_i=A\cap
\oplus_{|\alpha |\leq i}(\S1 A)^\d x^\alpha$ for all $i\geq 0$.

Fix elements $a_1, \ldots , a_s\in S$, then the derivations
$\d_1':= a_1\d_1, \ldots , \d_s':= a_1\d_s$ of $A$ are commuting
and locally nilpotent with the corresponding  filtration $\{
N_i'\}$ on $A$. Then $N_i'=N_i$ for all $i\geq 0$ where $\{ N_i\}$
is the filtration on $A$ determined by the derivations $\d_1,
\ldots , \d_s$.

More generally, fix elements $a_1, \ldots , a_s, t_1, \ldots ,
t_s\in S$, then consider  derivations $\d_1':= t_1^{-1}a_1\d_1,
\ldots , \d_s':= t_s^{-1}a_1\d_s$ of $\S1 A$ which are obviously
commuting and locally nilpotent and $\d_i'(A')\subseteq A'$ for
all $i$ where $A':= A[t_1^{-1}, \ldots , t_s^{-1}]$. Let
$A'=\cup_{i\geq 0}N_i'$ be the corresponding filtration associated
with the derivations $\d_1'\ldots, \d_s'$. Then, by Lemma
\ref{r22Dec05}, $N_i'=A'\cap \oplus_{|\alpha | \leq i} (\S1 A)^\d
x^\alpha $ for all $i\geq 0$. Instead of $A'$ one can take any
$\d'$-invariant subalgebra of $\S1 A$.


\section{Generators and defining relations for ring of invariants of commuting automorphisms}

Let $A$ be an algebra over a field $K$, $\s \in \Aut_K(A)$, and
$\d $ be a $\s$-{\em derivation} of the algebra $A$: $\d (ab)=\d
(a)\s (b)+a\d (b)$ for all $a,b\in A$. We will assume that $\d \s
=\s \d$. Then an induction on $n$  yields 
\begin{equation}\label{sdernab}
\d^n(ab)=\sum_{i=0}^n{n\choose i}\d^i(a)\s^i\d^{n-i}(b), \;\;
n\geq 1.
\end{equation}
It follows that the $A^\d := \ker (\d )$ is a subalgebra (of
constants for $\d $) of  $A$, or the ring of invariants, the union
$M:= M(\d , A)=\cup_{i\geq 0}M_i$ of the vector spaces $M_i:= \ker
(\d^{i+1})$ is a positively filtered algebra ($M_iM_j\subseteq
M_{i+j}$ for all $i,j\geq 0$), $M_0=A^\d \subseteq M_1\subseteq
\cdots $. For each $0\neq a\in M$, there exists a unique natural
number, say $d$, such that $a\in M_d\backslash M_{d-1}$. The
$d:=\deg (a)=\deg_\d (a)$ is called the $\d$-{\em degree} of the
element $a$.

{\it Example}. Given $\s\in \Aut_K(A)$, then $\d := \s -1$ is a
$\s $-derivation of the algebra $A$ such that $\d\s =\s\d$.

Given a vector space $V$ over the field $K$, a $K$-linear map $\v
: V\ra V$ is called {\em locally nilpotent} if, for all $v\in V$,
$\v^n(v)=0$ for all $n\gg 1$.

Given {\em commuting} $K$-automorphisms $\s_1,\ldots , \s_s$ of
the algebra $A$ such that the maps $\s_1-{\rm id}_A,\ldots ,
\s_s-{\rm id}_A$ are {\em locally nilpotent}. Then the maps
$\s_1-{\rm id}_A,\ldots , \s_s-{\rm id}_A$ are {\em commuting
locally nilpotent} $\s_1-, \ldots , \s_s-$derivations
respectively, and all the maps $\s_i$, $\s_j-{\rm id}_A$ {\em
commute}. The algebra $A$ has the filtration $\{ M_i\}_{i\geq 0}$
where $M_i:= \{ a\in A\, | \, (\s -{\rm id}_A)^\alpha (a)=0$ for
all $\alpha \in \Ns$ such that $|\alpha |>i\}$ where $(\s -{\rm
id}_A)^\alpha:= \prod_{i=1}^s (\s_i-{\rm id}_A)^{\alpha_i}$.
Clearly, $M_0=A^\s := \{ a\in A \, | \, \s_1(a)=\cdots =
\s_s(a)=a\}$, the {\em ring of} $\s$-{\em invariants},
$M_0\subseteq M_1\subseteq \cdots \subseteq M_i\subseteq \cdots
\subseteq  A=\cup_{i\geq 0}M_i$, $M_iM_j\subseteq M_{i+j}$ for all
$i,j\geq 0$ (use (\ref{sdernab})).

{\it Example}. Let $\s_1, \ldots , \s_s$ be $K$-automorphisms of
the polynomial algebra $P_s:= K[x_1, \ldots ,$ $ x_s]$ given by
the rule $\s_i(x_j)=x_j+\d_{ij}$. The automorphisms $\s_i$ commute
and all the maps $\s_i- {\rm id}_{P_s}$ are locally nilpotent.
Then the filtration $\{ M_i\}$ on the polynomial algebra $P_s$ is
the ordinary filtration: $M_i=\sum_{|\alpha |\leq i}Kx^\alpha$.

\begin{lemma}\label{c18Dec05}
Let $A$, $\d_1, \ldots , \d_s$, and $x_1, \ldots , x_s$  be as in
Theorem \ref{18Dec05}, and $\s \in \Aut_K(A)$. Then the
automorphism commute with the derivations $\d_1, \ldots , \d_s$
iff $\s (A^\d)=A^\d$ and $\s (x_i)=x_i+\l_i$ for some $\l_i\in
A^\d$.
\end{lemma}

{\it Proof}. $(\Rightarrow )$ If the automorphism $\s $ commutes
with derivations $\d_i$ then so does its inverse $\s^{-1}$, and so
$\s^{\pm 1} (A^\d )\subseteq A^\d$, hence $\s (A^\d )=A^\d$. By
(\ref{AAds}), $\s (x_i)=\l_i+\sum_{0\neq \alpha \in
\mathbb{N}^s}\l_{i, \alpha} x^\alpha$ for some $\l_i,\l_{i,
\alpha}\in A^\d$. Comparing the coefficients of $x^\alpha$'s in
the system of equations $\s \d_i (x_j)=\d_i\s (x_j)$, $1\leq
i,j\leq s$, yields  $\s (x_i)=\l_i +x_i$ for all $i$.

$(\Leftarrow )$ This implication is obvious.  $\Box $

\begin{theorem}\label{20Dec05}
Let $A$ be an arbitrary  $K$-algebra, $\s_1, \ldots , \s_s\in
\Aut_K(A)$ be  automorphisms of the algebra $A$. The following
statements are equivalent.
\begin{enumerate}
\item The maps  $\s_1-{\rm id}_A, \ldots , \s_s-{\rm id}_A$ are
commuting  locally nilpotent  and there exist elements $x_1,\ldots
, x_s\in A$ satisfying $\s_i(x_j)=x_j+\d_{ij}$ (the Kronecker
delta) for $1\leq i,j\leq s$. \item $\s_1=e^{\d_1}, \ldots ,
\s_s=e^{\d_s}$ for some commuting locally nilpotent derivations
$\d_1, \ldots , \d_s\in \Der_K(A)$ such that $\d_i(x_j)=\d_{ij}$,
$1\leq i,j\leq s$, for some elements $x_1,\ldots , x_s\in A$.
\end{enumerate}
If one of the two equivalent conditions holds then $\d_i:=
\sum_{k\geq 1}(-1)^{k+1}\frac{(\s_i-{\rm id}_A)^k}{k}$,  $A^\d
=A^\s : = \{ a\in \, | \, \s_1 (a)=\cdots = \s_s(a)=a\}$, and two
sets of $x$'s coincide up to adding elements of $A^\s =A^\d$. So,
one can apply all the previous results in finding generators and
defining relations for the algebra $A^\s$ (we leave it to the the
interested reader to write down the corresponding statements).
\end{theorem}

{\it Proof}. It is well-known that if an automorphism $\s \in
\Aut_K(A)$ is such that the map $\s -{\rm id}_A$ is a locally
nilpotent map then $\s =e^\d =\sum_{k\geq 0}\frac{\d^k}{k!}$ for a
unique locally nilpotent derivation $\d =\sum_{k\geq
1}(-1)^{k+1}\frac{(\s-{\rm id}_A)^k}{k}\in \Der_K(A)$, and vice
versa. It follows that the maps $\s_i-{\rm id}_A$ are commuting
locally nilpotent iff the derivations $\d_i$ are commuting locally
nilpotent. Then, $\s_i(x_j)=x_j+\d_{ij}$ iff $\d_i(x_j)=\d_{ij}$.
It is obvious that $A^\d =A^\s$.  $\Box$

{\it Example}. Let $P_n=K[x_1, \ldots , x_n]$ be a polynomial
algebra and $\s \in \Aut_K(P_n)$ where $\s (x_i)=x_i+x_{i-1}$,
$i=1, \ldots , n$, $x_0:=1$ (and $x_i:=0$ for all $i<0$). Then $\s
-{\rm id}$ is a locally nilpotent map and $P_n^\s =P_n^\d$ where
$\d := \sum_{k\geq 1}(-1)^k\frac{(\s -{\rm id})^k}{k}$, $\d
(x_1)=1$. By Theorem \ref{20Dec05} and Corollary \ref{ca18Dec05},
the ring of invariants $P_n^\s =K[\phi (x_2), \ldots , \phi
(x_n)]$ is a polynomial algebra in $n-1$ variables
\begin{eqnarray*}
 \phi (x_i)&=&\sum_{k\geq 0}(-1)^k\frac{x_1^k}{k!}\sum_{i_1, \ldots
 , i_k\geq 1} (-1)^{i_1+\cdots +i_k+k}\frac{(\s -{\rm id})^{i_1+\cdots +i_k}}{i_1\cdots i_k}
 \\
 &=& \sum_{k=0}^i\frac{x_1^k}{k!}\sum_{i_1, \ldots
 , i_k\geq 1} (-1)^{i_1+\cdots +i_k}\frac{x_{i-i_1-\cdots -i_k }}{i_1\cdots
 i_k}.
\end{eqnarray*}

\begin{corollary}\label{d25Dec05}
Let $A$ be an arbitrary algebra over the  field $K$. The following
statements are equivalent.
\begin{enumerate}
\item There exist commuting  $K$-automorphisms $\s_1,\ldots ,\s_s$
of the algebra $A$ such that the maps  $\s_i -{\rm id}$ are
locally nilpotent and $\s_i(x_j)=x_j+\d_{ij}$ for some elements
$x_1,\ldots , x_s\in A$.
 \item The algebra $A$ is an iterated Ore extension
 $A=B[x_1;d_1]\cdots [x_s;d_s]$ such that $d_i(B)\subseteq B$ and
 $d_i(x_j)\in B$ for all $1\leq i,j\leq s$.
\end{enumerate}
If, say, the first condition holds then $A=A^\s  [x_1;d_1]\cdots
[x_s;d_s]$  is an iterated Ore extension of the ring $A^\s$ such
that $d_i =\ad (x_i)$, $[x_i, A^\s ]\subseteq A^\s$, and
$[x_i,x_j]\in A^\s$ for all $i,j$.   In particular,
$A=\oplus_{\alpha \in \Ns}x^\alpha A^\s = \oplus_{\alpha \in \Ns}
A^\s x^\alpha=\cup_{i\geq 0}M_i$ where $ M_i=\oplus_{|\alpha |\leq
i} A^\s x^\alpha $.
\end{corollary}

{\it Proof}. $(1\Rightarrow 2)$ By Theorem \ref{20Dec05}, we have
a set $\d_1, \ldots , \d_s$ of commuting locally nilpotent
derivations of the algebra $A$ such that $\d_i(x_j)=\d_{ij}$ for
all $i,j$. By Theorem \ref{18Dec05}, statement 2 holds.

$(2\Rightarrow 1)$ Given the iterated Ore extension as in
statement 2. It is easy to check that the $K$-automorphisms $\s_1,
\ldots , \s_s\in \Aut_K(A)$ given by the rule
$\s_i(x_j)=x_j+\d_{ij}$ satisfy the conditions of statement 1. The
rest follows from Theorem \ref{18Dec05} and the fact that $(\s
-{\rm id})^\alpha (x^\beta )=\alpha ! \d_{\alpha , \beta}$ for all
$\alpha , \beta \in \Ns $ such that $|\alpha |\geq |\beta |$ where
$(\s -{\rm id})^\alpha := \prod_{i=1}^s(\s_i-{\rm id
})^{\alpha_i}$. $\Box $

Corollary \ref{d25Dec05} proves that the filtration $\{M_i\}$ of
the algebra $A$ for the automorphisms $\s_1, \ldots , \s_s$
coincides with the filtration $\{ N_i\}$ for the derivations
$\d_1, \ldots , \d_s$ (where $\s_i=e^{\d_i}$), that is $M_i=N_i$
for all $i\geq 0$.


\section{The inverse map for  automorphism that preserve
 the ring of invariants of derivations}
Let $A$, $\d_1, \ldots , \d_s$ and $x_1, \ldots , x_s$ be as in
Theorem \ref{18Dec05}. Suppose that an automorphism $\s\in
\Aut_K(A)$  preserves the ring of invariants $A^\d$ ($\s (A^\d
)=A^\d$). Let $\s_\d := \s|_{A^\d}\in \Aut_K(A^\d)$. Suppose we
know the inverse $\s_\d^{-1}$ and the twisted derivations $\d_1':=
\s \d_1\s^{-1}, \ldots , \d_s':= \s \d_s\s^{-1}\in \Der_K(A)$,
then we can write {\em explicitly} a formula for the inverse
automorphism $\s^{-1}$ of $\s $ (Theorem \ref{i19Dec05} and
Theorem \ref{i8Nov05}). Since $A=\oplus_{\alpha \in \Ns}A^\d
x^\alpha = \oplus_{\alpha \in \Ns} x^\alpha A^\d $, the
automorphism $\s$ is uniquely determined by its restriction
$\s_\d$ to the ring of invariants $A^\d$ and the images of $x$'s:
\begin{equation}\label{gxsh1}
x_1':= \s (x_1), \ldots , x_s':=\s (x_s).
\end{equation}
The twisted derivations $\d_1':= \s \d_1\s^{-1}, \ldots , \d_s':=
\s \d_s\s^{-1}\in \Der_K(A)$ is a set of {\em commuting locally
nilpotent} derivations of the algebra $A$ satisfying
$\d_i'(x_j')=\d_{ij}$. For each $i=1, \ldots , s$, consider the
maps
$$ \phi_i':=\sum_{k\geq 0}(-1)^k\frac{(x_i')^k}{k!}(\d_i')^k,  \;\; \psi_i':=\sum_{k\geq 0}(-1)^k(\d_i')^k
(\cdot ) \frac{(x_i')^k}{k!}\, : \, A\ra A$$ which are
homomorphisms of right and left $A^\d$-modules respectively.
 The maps
\begin{eqnarray*}
\phi_\s &: =&\phi_s'\phi_{s-1}'\cdots \phi_1':A\ra A, \;\;
a=\sum_{\alpha \in \mathbb{N}^s}  x^\alpha_\alpha a_\alpha \mapsto
\phi_\s
(a)=a_0,\\
\psi_\s &:=& \psi_1'\psi_2'\cdots \psi_s':A\ra A, \;\;
a=\sum_{\alpha \in \mathbb{N}^s} a_\alpha x^\alpha\mapsto \psi_\s
(a)=a_0,
\end{eqnarray*}
are {\em projections} onto the subalgebra $A^\d$ of $A=A^\d \oplus
(\oplus_{0\neq \alpha \in \Ns }x^\alpha A^\d )$ and $A=A^\d \oplus
(\oplus_{0\neq \alpha \in \Ns } A^\d x^\alpha )$ respectively,
they are homomorphisms of right and left $A^\d $-modules
respectively.
 By Theorem \ref{a18Dec05}, for any $a\in A$,
$$ a=\sum_{\alpha \in \mathbb{N}^s}(x')^\alpha \phi_\s  (\frac{(\d')^\alpha}{\alpha
!} a)=\sum_{\alpha \in \mathbb{N}^s}\psi_\s
(\frac{(\d')^\alpha}{\alpha !} a)x^\alpha.$$ Then applying
$\s^{-1}$ we finish the proof of the next theorem.

\begin{theorem}\label{i19Dec05}
Let $A$, $\d_i$, $\d_i'$, $x_i$, and $x_i'$ be as above (i.e. the
algebra $A$ is from Theorem \ref{18Dec05}). For $a\in A$,
$$ \s^{-1}(a)= \sum_{\alpha \in \Ns} x^\alpha \s_\d^{-1}\phi_\s (\frac{(\der')^\alpha}{\alpha
!}a)= \sum_{\alpha \in \Ns}  \s_\d^{-1}\psi_\s
(\frac{(\der')^\alpha}{\alpha !}a)x^\alpha.$$
\end{theorem}

As an application of Theorem \ref{i19Dec05} we find the inverse
map of an automorphism of the Weyl algebra with polynomial
coefficients.

 The {\em Weyl } algebra $A_n=A_n(K)$ is a
$K$-algebra generated by $2n$ generators $x_1, \ldots , x_{2n}$
subject to the defining relations:
$$ [x_{n+i}, x_j]=\d_{ij}, \;\; [x_i, x_j]=[x_{n+i}, x_{n+j}]=0\;\; {\rm
for\;\; all}\;\; 1\leq i,j\leq n,$$ where $\d_{ij}$ is the
Kronecker delta, $[a,b]:=ab-ba$.

Let $P_m=K[x_{2n+1}, \ldots x_{2n+m}]$ be  a polynomial algebra.
The Weyl algebra with polynomial coefficients $A:= A_n\t
P_m=\bigoplus_{\alpha \in \mathbb{N}^s} Kx^\alpha $ where
$s:=2n+m$, $x^{\alpha }:= x_1^{\alpha_1 }\cdots x_s^{\alpha_s }$,
the {\em order} of the $x$'s in the product  is {\em fixed}. The
algebra $A_n\t P_m$ admits the finite set of {\em commuting
locally nilpotent} derivations, namely, the `partial derivatives':
$ \der_1:= \frac{\der }{\der x_1}, \ldots , \der_s:= \frac{\der
}{\der x_s}$.
 Clearly, $\der_i=\ad (x_{n+i})$ and $\der_{n+i}=-\ad (x_i)$,
$i=1,\ldots , n$.

Let  $\Aut_K(A_n\t P_m)$ be the group of $K$-algebra automorphisms
of the algebra $A_n\t P_m$. Given an automorphism $\s
\in\Aut_K(A_n\t P_m)$. It is uniquely determined by the elements $
x_1':= \s (x_1), \ldots , x_s':=\s (x_s) $ of the algebra $A_n\t
P_m$. The centre $Z:=Z(A_n\t P_m)$ of the algebra $A_n\t P_m$ is
equal to $P_m$. Therefore, the restriction $\s|_{P_m}\in
\Aut_K(P_m)$, and so $ \D :=\det (\frac{\der x_{2n+i}'}{\der
x_{2n+j}})\in K^*$  where $i,j=1, \ldots , n$. The corresponding
(to the elements $x_1',\ldots , x_s'$) `partial derivatives' (the
set of commuting locally nilpotent derivations of the algebra
$A_n\t P_m$) 
\begin{equation}\label{xsh2}
\der_1':= \frac{\der }{\der x_1'}, \ldots , \der_s':= \frac{\der
}{\der x_s'}
\end{equation}
are equal to 
\begin{equation}\label{dad1}
\der_i':= \ad (\s (x_{n+i})), \;\; \der_{n+i}':= -\ad (\s
(x_{i})), \;\; i=1, \ldots , n,
\end{equation}

\begin{equation}\label{dad2}
\der_{2n+j}'  := \D^{-1} \det
 \begin{pmatrix}
  \frac{\der \s (x_{2n+1})}{\der x_{2n+1}} & \cdots & \frac{\der \s (x_{2n+1})}{\der x_{2n+m}} \\
  \vdots & \vdots & \vdots \\
\frac{\der }{\der x_{2n+1}} & \cdots & \frac{\der }{\der x_{2n+m}}\\
 \vdots & \vdots & \vdots \\
\frac{\der \s (x_{2n+m})}{\der x_{2n+1}} & \cdots & \frac{\der \s (x_{2n+m})}{\der x_{2n+m}} \\
\end{pmatrix}, \;\;\; j=1, \ldots , m,
\end{equation}
where we `drop' $\s (x_{2n+j})$ in the determinant $\det
(\frac{\der \s (x_{2n+k})}{\der x_{2n+l}})$. Clearly, $\der_i'=\s
\der_i\s^{-1}$ for $i=1, \ldots , s$, and $A^\der =K$, so
$\s|_K={\rm id}_K$ is known. Now, one can apply Theorem
\ref{i19Dec05} to have the next result.

\begin{theorem}\label{i8Nov05}
\cite{inform'05} {\rm (The Inversion Formula)} For each $\s \in
\Aut_K(A_n\t P_m)$
 and $a\in A_n\t P_m$,
 $$ \s^{-1}(a)=\sum_{\alpha \in \mathbb{N}^s}x^\alpha \phi_\s
 (\frac{(\der')^\alpha}{\alpha!}a)=\sum_{\alpha \in \mathbb{N}^s}\psi_\s
 (\frac{(\der')^\alpha}{\alpha!}a)x^\alpha , $$
 where $(\der')^{\alpha} :=(\der_1')^{\alpha_1}\cdots
 (\der_s')^{\alpha_s}$ and  $s=2n+m$.
\end{theorem}


\section{Integral closure and commuting locally nilpotent derivations}

In this section, the structure of algebras is described that admit
a set of commuting locally nilpotent derivations with left
localizable  kernels.

For an arbitrary algebra $A$, we say that derivations $\d_1,\ldots
, \d_s$ of the algebra $A$ have {\em generic kernels} iff the
$2^s-1$ sets $\{ \cap_{i\in I} A^{\d_i}\, | \, \emptyset \neq
I\subseteq \{ 1, \ldots , s\}\}$  are {\em distinct} (iff the sets
$A^\d$, $\Adhi := \cap_{j\neq i}A^{\d_j}$, $i=1,\ldots , s$ are
distinct iff $A^\d \neq \Adhi $ for $i=1,\ldots , s$). We say that
the derivations $\d_1,\ldots , \d_s$ of the algebra $A$ have {\em
left localizable kernels} iff there exists a {\em left Ore} set
$S$ of the algebra $A$ such that $S\subseteq A^\d \cap A^{reg}$
and $S\cap \d_i(\Adhi )\neq \emptyset$ for all $i=1,\ldots ,s$,
where $A^{reg}$ is the set of all {\em regular} elements of the
algebra $A$ (an element $a\in A$ is {\em regular} if, by
definition, it is not a left and right zero divisor of the algebra
$A$). {\em If the derivations $\d_1, \ldots , \d_s$ have left
localizable kernels then they have generic kernels}: for each $i$,
fix $y_i\in \Adhi$ such that $\d_i (y_i)\in S$, then $y_i\in \Adhi
\backslash (A^\d \cup \cup_{j\neq i} \Adhj )$, and so the
derivations have generic kernels.  Clearly, if there exists
elements $x_1, \ldots , x_s\in A$ such that $\d_i(x_j)=\d_{ij}$
then the derivations $\d_1, \ldots , \d_s$ have left localizable
kernels (but not vice versa): for take $S=\{ 1\}$.

 \begin{theorem}\label{24Dec05}
Let $A$ be an (arbitrary) algebra over the field $K$. The
following statements are equivalent.
\begin{enumerate}
\item The algebra $A$ admits a finite set of commuting locally
nilpotent derivations, say $\d_1, \ldots , \d_s\in \Der_K(A)$,
with left localizable  kernels. \item There exists a left Ore set
$S$ of the algebra $A$ such that $S\subseteq A^{reg}$, $\S1 A =
B[x_1; d_1]\cdots $ $  [ x_s; d_s]$ is an iterated Ore extension
such that $S\subseteq B$, $d_i (B)\subseteq B$ and $d_i(x_j)\in B$
for all $1\leq i,j\leq s$, and the algebra $A$ is
$\der_i$-invariant $(\der_i(A)\subseteq A)$ for all $1\leq i \leq
s$ where $\der_i := \frac{\der}{\der x_i}\in \Der_B(\S1 A)$ are
the formal partial derivatives of the $B$-algebra  $\S1 A$.
\end{enumerate}
If, say, the first condition holds then there exists a left Ore
set $S\subseteq A^\d \cap A^{reg}$ such that $\S1 A=(\S1 A)^\d
[x_1; d_1]\cdots [ x_s; d_s]$ is an iterated Ore extensions such
that $d_i ((\S1 A)^\d )\subseteq (\S1 A)^\d $ and $d_i(x_j)\in
(\S1 A)^\d $ for all $1\leq i,j\leq s$. In particular, $\S1 A=
\oplus_{\alpha \in \Ns } (\S1 A)^\d x^\alpha = \oplus_{\alpha \in
\Ns } x^\alpha (\S1 A)^\d $ and $\S1 A=\cup_{i\geq0}N_i'$,
$N_i'=\oplus_{|\alpha |\leq i } (\S1 A)^\d
x^\alpha=\oplus_{|\alpha |\leq i } x^\alpha (\S1 A)^\d $, $i\geq
0$. Finally, $A=\cup_{i\geq 0}N_i$ and $N_i'=\S1 N_i$ for all
$i\geq 0$, in particular, $\S1 (A^\d)=(\S1 A)^\d $ and $A^\d
=A\cap (\S1 A)^\d$.
\end{theorem}

{\it Proof}. $(1\Rightarrow 2)$  The derivations $\d_i$ are left
localizable, that is, there exists a left Ore set $S$ of the
algebra $A$ such that $S\subseteq A^\d \cap A^{reg}$ and $S\cap
\d_i (\Adhi )\neq \emptyset$ for all $i=1, \ldots , s$. So, for
each $i=1, \ldots , s$, one  can pick up an element, say $y_i\in
\Adhi$,  such that $\d_i(y_i)\in S$, then for the elements $x_i:=
\d_i(y_i)^{-1}y_i\in \S1 A$ we have $\d_j(x_i)=\d_{ij}$, where the
`new' derivation $\d_j$  is the unique extension of the `old'
derivation $\d_j$ to a derivation of the algebra $\S1 A$. By
Theorem \ref{18Dec05}, $\S1 A= B[x_1;d_1]\cdots [ x_s; d_s]$ is an
iterated Ore extension such that $ B=(\S1 A)^\d$, $d_i(B)\subseteq
B$,  $d_i(x_j)\in B$ for all $i,j$, and $\d_i:= \frac{\der}{\der
x_i}\in \Der_B(\S1 A)$ are formal partial derivatives over $B$.
Now, it is obvious that the algebra $A$ is $\frac{\der}{\der
x_i}$-invariant for all $i$.

$(2\Rightarrow 1)$ Suppose that the second statement holds. The
derivations $\der_1, \ldots , \der_s\in \Der_B(\S1 A)$ are
commuting locally nilpotent, hence so are their restrictions, say
$\d_1, \ldots , \d_s$, to the $\der$-invariant subalgebra $A$ of
$\S1 A$ (the set $S$ consists of regular elements of the algebra
$A$, so one can identify the algebra $A$ with its isomorphic image
in $\S1 A$  under the natural monomorphism $A\ra \S1 A$, $a\mapsto
\frac{a}{1}$). For each $x_i$, fix an element $s_i\in S$ such that
$y_i:= s_ix_i\in A$. Then $\d_j (y_i)=s_i\d_{ij}$ for all $i,j$.
Since $S\subseteq A^{reg}\cap A^\d$ is a left Ore set and $\d_i
(y_i)=s_i\in \d_i (\Adhi )\cap S$, the kernels of the derivations
$\d_1, \ldots , \d_s$ are left localizable. This finishes the
proof of the implication.

The rest is a direct consequence of Theorem \ref{18Dec05} and the
fact that $S\subseteq A^{reg}\cap A^\d$.  $\Box $

 \begin{corollary}\label{c24Dec05}
Let a $K$-algebra  $A$ be a commutative domain. The following
statements are equivalent.
\begin{enumerate}
\item The algebra $A$ admits a finite set of commuting locally
nilpotent derivations, say $\d_1, \ldots , \d_s\in \Der_K(A)$,
with generic kernels. \item There exists a nonzero element $t\in
A$ such that the localization $A_t:= A[t^{-1}]$ of the algebra $A$
at the powers of the element $t$ is a polynomial algebra
$A_t=B[x_1, \ldots , x_s]$ such that $t\in B$ and the algebra $A$
is $\der_i$-invariant ($\der_i(A)\subseteq A$) for all $1\leq i
\leq s$ where $\der_i := \frac{\der}{\der x_i}\in \Der_B(A_t)$ are
the formal partial derivatives of  $ A_t$ over $B$.
\end{enumerate}
\end{corollary}

{\it Proof}. $(1\Rightarrow 2)$ The derivations $\d_i$ have
generic kernels, so the algebras $A^\d$ and $\Adhi$, $i=1,\ldots ,
s$, are distinct. So, for each $i=1, \ldots , s$, one  can fix an
element, say $y_i\in \Adhi$, such that $0\neq \d_i(y_i)\in A^\d$.
 Then the element $t:= \d_1(y_1)\cdots \d_s(y_s)\in A^\d$ is a
 nonzero one since the algebra $A$ is a domain, and the
 derivations $\d_1, \ldots , \d_s$ have (left) localizable generic
 kernels, for it suffices to take $S=\{ t^i\, | \, i\geq 0\}$.
 Applying Theorem \ref{24Dec05}, we obtain statement 2.

$(2\Rightarrow 1)$ This implication is obvious because of  Theorem
\ref{24Dec05}. $\Box $

The next result gives explicitly generators and defining relations
for the integral closure $\tK$ of the field $K$ in the algebra
$A$.
 \begin{corollary}\label{1c24Dec05}
Let a domain $A=K\langle y_1,\ldots y _r\rangle$ be an affine
commutative $K$-algebra of Krull dimension $s\geq 1$, $\tK$ be an
algebraic closure of the field $K$ in the algebra $A$ ($\tK $ is a
field finite over $K$, i.e. $[\tK :K]<\infty$). The following
statements are equivalent.
\begin{enumerate}
\item There exist $s$ commuting locally nilpotent derivations, say
$\d_1, \ldots , \d_s\in \Der_K(A)$, with generic kernels. \item
$A=\tK [x_1,\ldots , x_s]$ is a polynomial algebra over the field
$\tK$ in $s$ variables.
 \item There exist derivations $\d_1, \ldots , \d_s\in \Der_{\tK
 }(A)$ and elements $x_1, \ldots , x_s\in A$ such that
 $\d_i(x_j)=\d_{ij}$ (the Kronecker delta) for all $1\leq i,j\leq
 s$.
 \end{enumerate}
If, say, the first statement holds and $\{ N_i\}$ be the
filtration on $A$ associated with the derivations $\d_1,\ldots ,
\d_s$ then

$(i)$ $\dim_K(N_i)=[\tK :K]{i+s\choose s}= \frac{[\tK
:K]}{s!}i^s+\cdots $ for all $i\geq 0$.

$(ii)$ The map $\phi : A\ra \tK$ (from Corollary \ref{ca18Dec05})
is an algebra epimorphism with kernel $(x_1, \ldots , x_s)$, i.e.
$\tK \simeq A/(x_1,\ldots , x_s)$. Alternatively, the field $\tK$
is generated over $K$ by the elements $\phi (y_1), \ldots , \phi
(y_r)$ that satisfy the defining relations $\CR = \{ f(t_1, \ldots
, t_r)\in K[t_1, \ldots , t_r]\, | \, f(\phi (y_1), \ldots , \phi
(y_r))\in (x_1,\ldots , x_s)\}$.
\end{corollary}

{\it Proof}. $(1\Rightarrow 3)$ By Corollary \ref{c24Dec05}, $A_t
=A_t^\d [ x_1, \ldots , x_s]$ for a nonzero element $t\in A^\d$.
 Since $s=\Kdim (A)=\Kdim(A_t)$, we see that $A_t^\d$ is a field
 (since $A$ is a domain) which is finite over $K$, i.e.
 $[A^\d_t:K]<\infty$. The element $t\in A^\d \subseteq A^\d_t$ is
 algebraic, hence $t^{-1}\in A^\d$, and so $A_t=A$ and $A=A^\d
 [x_1, \ldots , x_s]$. Clearly, $A^\d \subseteq \tK$. The
 reverse inclusion is obvious  since char$(K)=0$ (if $u\in \tK$
 then $f(u)=0$, $f'(u):=\frac{df}{dx}(u)\neq 0$ for some nonzero
 polynomial $f(x)\in K[x]$, and so $0=\d_i (f(u))=f'(u)\d_i(u)$
 implies $\d_i(u)=0$ for all $i$ which means that $u\in A^\d$),
 hence $A^\d =\tK $. It is obvious that $\d_i=\frac{\der}{\der
 x_i}\in \Der_{\tK }(A)$ and $\d_i(x_j)=\d_{ij}$.

$(3\Rightarrow 2)$ By Theorem \ref{18Dec05}, $A=A^\d [x_1,\ldots ,
x_s]$ is a polynomial algebra with coefficients from the algebra
of invariants $A^\d$. Repeating the above argument we have $A^\d
=\tK$.

$(2\Rightarrow 1)$ The formal partial derivatives
$\d_i:=\frac{\der}{\der x_i}\in \Der_{\tK }(A)$ are commuting
locally nilpotent derivations with generic kernels since $x_i\in
\Adhi \backslash (A\cup \cup_{j\neq i} \Adhj )$. This finishes the
proof of the equivalence of the three statements.

Suppose that the equivalent conditions hold, then
$N_i=\oplus_{|\alpha |\leq i} \tK x^\alpha$, and so
 $\dim_K(N_i)=[\tK :K]{i+s\choose s}= \frac{[\tK
:K]}{s!}i^s+\cdots $ for all $i\geq 0$, i.e. $(i)$ is proved. The
statement $(ii)$ follows from Corollary   \ref{ca18Dec05}. $\Box $

{\it Remark}. The finite separable field extension $\tK /K$ is
generated by a single element, say $x$, over $K$. So, the algebra
$A$ from Corollary \ref{1c24Dec05} is generated by $s+1$ elements
$x, x_1,\ldots , x_s$ that satisfy a single defining relations
$f(x)=0$ for an irreducible polynomial $f(y)\in K[y]$ of degree
$[\tK :K]$.


\section{A construction of simple algebras}

In this section, a construction of simple algebras is given
(Theorem \ref{26Dec05}) that comes from a set of commuting locally
nilpotent derivations which satisfy the conditions of Theorem
\ref{18Dec05}.

\begin{theorem}\label{26Dec05}
Let $A$, $\d_1, \ldots , \d_s$ and $x_1, \ldots , x_s$ be as in
Theorem \ref{18Dec05}. Given a (two-sided) maximal  ideal $\gm$ of
the algebra $A^\d$ such that $[x_i, \gm ]\subseteq \gm $ and
$[x_i, Z]\subseteq \gm$ for all $i=1, \ldots , s$ where $Z$ is the
centre of the factor algebra $A^\d / \gm$. Then the iterated Ore
extension $\CA := A/(\gm )[t_1, \ldots , t_s; \d_1, \ldots ,
\d_s]$ of the algebra $A/(\gm )$  is a simple algebra where $(\gm
):= A\gm A$, the elements $t_1, \ldots , t_s$ commute, and
$t_ia=at_i+\d_i(a)$ for all $a\in A/(\gm )$ where $\d_i\in
\Der_K(A/(\gm ))$ is the induced derivation: $u+(\gm )\mapsto
\d_i(u)+(\gm )$, $u\in A$.
\end{theorem}

{\it Proof}. Using Theorem \ref{18Dec05} and abusing notation
slightly  one can write the factor algebra $A/(\gm )$ as the
iterated Ore extension $A^\d /\gm [x_1;d_1]\cdots [ x_s; d_s]$
 of the algebra $A^\d /\gm$. So, without loss of generality we can assume that $\gm
=0$, that is $A^\d$ is a simple algebra. We have to prove that the
iterated Ore extension $\CA :=A [t_1, \ldots , t_s: \d_1, \ldots ,
\d_s]$ of the algebra $A$ is a simple algebra. The algebra $A^\d$
is simple,and so its centre $Z$ is a field that contains the field
$K$. Let $I$ be a nonzero ideal of the algebra $\CA$, we have to
show that $I=\CA$. Recall that $\CA =\oplus_{\alpha \in
\Ns}At^\alpha$ where $t^\alpha:= t_1^{\alpha_1}\cdots
t_s^{\alpha_s}$ and $ A=\oplus_{\alpha \in \Ns} A^\d x^\alpha$.
Fix a nonzero element, say $a\in I$. Then $a=\sum a_\alpha
t^\alpha$ for some elements $a_\alpha \in A$ not all of which are
zero. Note that the inner derivation $\ad (t_i)$ of the algebra
$\CA$ is a formal partial derivative $\frac{\der}{\der x_i}$ over
$A^\d$ of the algebra $\CA = A^\d \langle x_1, \ldots , x_s, t_1,
\ldots , t_s\rangle$, that is $\frac{\der}{\der x_i}(A^\d )=0$,
$\frac{\der}{\der x_i}(x_i)=1$ and $\frac{\der}{\der x_i}(y)=0$
for all $y\in \{ x_1, \ldots , \widehat{x}_i, \ldots , x_s, t_1,
\ldots , t_s\}$ (the hat over a symbol means that it is missed).
Note that the ideal $I$ is $\frac{\der}{\der x_i}$-invariant for
all $i=1, \ldots , s$. Applying carefully several times inner
derivations of the type $\ad (t_i)=\frac{\der}{\der x_i}$ to the
element $a$ we see that we can assume that all the coefficients
$a_\alpha \in A^\d$ and not all of which are zero ones. Let
$V\subseteq A^\d$ be the vector space over the field $Z$ generated
by all the coefficients $a_\alpha$. Suppose that a set $a_\alpha ,
a_\beta , \ldots , a_\g$ is a $Z$-basis for $V$. By {\em the
Density Theorem}, there
 are elements $u_1, \ldots u_k, v_1, \ldots , v_k\in A^\d$ such that
 $\sum_{i=1}^k u_ia_\alpha v_i=1$, $\sum_{i=1}^k u_ia_\beta v_i=0, \ldots , \sum_{i=1}^k u_ia_\g
 v_i=0$. Applying the map $A \ra A$, $(\cdot )\mapsto
 \sum_{i=1}^k u_i(\cdot )v_i$, to the element $a$, we can assume
 that all the coefficients $a_\alpha \in Z$ but not all are zero.
  By the assumption,  $[x_i, Z]=0$ for all $i=1, \ldots , s$. Then applying
 carefully the inner derivations of the type $ -\ad (x_i)$ to the
 element $a$ and taking into account the fact that $-\ad
 (x_i)(t_j)=\d_{ij}$, we get an element $0\neq b\in Z\cap I$.
 Hence, $I=\CA$, as required.
$\Box $

{\it Example}. Let $A=P_n:=K[x_1, \ldots , x_n]$,
$\d_1:=\frac{\der}{\der x_1}, \ldots , \d_s:=\frac{\der}{\der
x_s}\in \Der_K(P_n)$, $A^\d =K$, $\gm =0$. Then the algebra
$P_n[t_1, \ldots , t_s; \d_1, \ldots , \d_s]$ is the $n$'th Weyl
algebra $A_n$.

{\it Example}. Let $A=F_2:=K\langle x_1, x_2\rangle$ is the free
algebra,  $\d_1:=\frac{\der}{\der x_1},  \d_2:=\frac{\der}{\der
x_2}\in \Der_K(F_2)$, the ideal  $\gm $ of $F_2^\d$ is generated
by a single element $[x_2, x_1]-1$ is $\ad (x_i)$-invariant, $i=1,
2$. Then $F_2^\d /\gm \simeq K$ and the algebra $\CA =K[x_1][x_2;
d_2:=\frac{\der}{\der x_1}][t_1, t_2; \d_1, \d_2]$ is a simple
algebra.

{\it Example}. Let $A:=F_s=K\langle x_1, \ldots x_s\rangle$,
$s\geq 2$, be a free algebra, $\d_1:=\frac{\der}{\der x_1},\ldots
, \d_s:=\frac{\der}{\der x_s}$. Let $I$ be an ideal of the algebra
$F_s^\d$ generated by all the commutators $[x_i, [x_j,x_k]]$. Then
the factor algebra  $P:=F_s^\d / I$ is a polynomial algebra in
${s\choose 2}$ variables $y_{ij}:= [x_i,x_j]+I$ and $\bA := A/(I)=
P[\bx_1 ; \ad \, \bx_1]\cdots [\bx_s ; \ad \, \bx_s] $ (see
Corollary \ref{f26Dec05}). Note that all $(\ad \, \bx_i)|_{P}=0$
and $[\bx_i, \bx_j ]=y_{ij}$. Hence, every maximal ideal $\gm $ of
the algebra $P$ satisfies the conditions of Theorem \ref{26Dec05}
and one can easily see that the algebra $\CA $ is isomorphic to
the Weyl algebra $A_s$ over the field $L:= P/\gm $. Note that the
factor algebra $\bA / (\gm )$ is isomorphic to the tensor product
$A_s\t_LP_m$ of the Weyl algebra $A_s$ and a polynomial algebra
$P_m$ in $m$ variables such that $s=2n+m$ and $2n$ is the rank of
the $s\times s$ skew  symmetric matrix  $(y_{ij}+\gm )$ over $L$.

{\it Example}. The same results are true for a {\em free
metabelian algebra}. Let $J$ be an ideal of the free algebra
$F_s$, $s\geq 2$, generated by all the double commutators $[a,
[b,c]]$ where $a,b,c\in F_s$. The ideal $J$ is $\d $-invariant.
Hence, (by definition) the {\em free metabelian algebra} is
defined as  $R:= F_s/J$ and it is isomorphic to the factor algebra
of $\bA$ (from the previous example) by an ideal $(J')$ generated
by an ideal $J'$ of the polynomial algebra $P$, i.e. $R\simeq P/J'
[\bx_1 ; \ad \, \bx_1]\cdots [\bx_s ; \ad \, \bx_s]$. Now, it is
obvious (it is a particular case of the previous example) that,
for any maximal ideal $\gm $ of the algebra $P/J'$, the algebra
$\CA$ is isomorphic to the Weyl algebra $A_s$ over the field $L:=
(P/J')/\gm $, and the factor algebra $R/(\gm )$ is isomorphic to
the tensor product $A_n\t_LP_m$, $s=2n+m$ as in the example above.


\section{Linear maps as differential operators}
Let $A:= A_n\t P_m=\oplus_{\alpha \in \mathbb{N}^s} Kx^\alpha$,
$s:=2n+m$, be the $n$'th Weyl algebra with polynomial coefficients
$P_m$. The set of formal `partial derivatives' $\der_1:=
\frac{\der}{\der x_1}, \ldots , \der_s:= \frac{\der}{\der x_s}$ is
a set of {\em commuting locally nilpotent $K$-derivations} of the
algebra $A$.  Consider the algebra $\hA := A[[\der_1, \ldots ,
\der_s]]=\oplus_{\alpha \in \mathbb{N}^s}A\der^\alpha$,
$\der^\alpha :=\der_1^{\alpha_1}\cdots \der_s^{\alpha_s}$, of
formal (noncommutative) series $\sum a_\alpha \der^\alpha$,
$a_\alpha \in A$, with multiplication given by the rule $\der_i a
=a\der_i + \der_i(a)$, $a\in A$, $1\leq i\leq s$. The
multiplication of series is well-defined since all the derivations
commute and are locally nilpotent. Since $\der_i \in
\Der_K(A)\subseteq {\rm End}_K(A)$, the algebra $\hA$ is, in fact,
a subalgebra of the algebra ${\rm End}_K(A)$ of all $K$-linear
endomorphisms of the vector space $A$. The next theorem shows that
they coincide.

\begin{theorem}\label{16Dec05}
$\hA ={\rm End}_K(A)$.
\end{theorem}

{\it Proof}. The algebra $A$ has a natural finite dimensional
filtration $\{ \CA_i := \sum_{|\alpha |\leq i}Kx^\alpha \}_{i\geq
0}$ $(\CA_i\CA_j\subseteq \CA_{i+j}$ for all $i,j\geq 0$), $\dim
(\CA_i)={s+i\choose s}$, and $\der^\alpha (\CA_i)\subseteq
\CA_{i-|\alpha |}$ for all $ \alpha \in \mathbb{N}^s$ and $i\geq
0$ (we set $\CA_i:=0$ for negative $i$). We have mentioned in
passing that the algebra $\hA$ is a subalgebra of ${\rm
End}_K(A)$, let us prove this statement, that is {\em each nonzero
series} $a=\sum a_\alpha \der^\alpha\in \hA$ {\em determines a
nonzero linear map}: let $i:= \min \{ | \alpha | \, | \, a_\alpha
\neq 0 \}$, fix $a_\alpha \neq 0$ with $|\alpha | = i$, then
$a(x^\alpha )=a_\alpha \alpha !\neq 0$, as required.

It remains to show that that any linear map $f\in {\rm End}_K(A)$
can be represented by a series $a=\sum a_\alpha \der^\alpha\in
\hA$. It means that $ f(x^\beta )= a(x^\beta )$, for  all $\beta
\in \mathbb{N}^s$.  The unknowns coefficients $a_\alpha \in A$ can
be found from this system step by step. Clearly, $f(1)=a_0$.
 Suppose that $i>0$ and all the coefficients $a_\alpha $ with
$|\alpha | <i$ have been found. Then, for each $\alpha$ such that
$|\alpha | =i$, the element $a_\alpha$ can be found (uniquely)
from the equation $f(x^\alpha )=\alpha ! a_\alpha +\sum_{|\beta |
<i}\der^\beta (x^\alpha )$.  $\Box $

Now we are ready to give a short direct proof of the fact that
$A_n=\CD (P_n)$. 
\begin{corollary}\label{An=DPn}
Let $K$ be a field of characteristic zero. The Weyl algebra $A_n$
is  the ring of differential operators $\CD (P_n)$ with polynomial
coefficients.
\end{corollary}

{\it Proof}. Applying Theorem \ref{16Dec05} to the polynomial
algebra $A=P_n=K[x_1, \ldots , x_n]$, we have ${\rm
End}_K(P_n)=\widehat{P}_n$. Let
$\mathbb{N}=\oplus_{i=1}^n\mathbb{N}e_i$ where $e_1=(1, 0, \ldots,
0), \ldots , e_n=(0,  \ldots, 0, 1)$. It follows from $[x_i,
\der^{\alpha}]=\alpha_i\der^{\alpha - e_i}$ for all $\alpha$ and
$i$ (and from definition of the ring of differential operators)
that the $j$'th term of the {\em order} filtration of the ring of
differential operators  $\CD (P_n)$ on $P_n$ is equal to
$\oplus_{|\alpha |\leq j}P_n\der^\alpha$. Hence $\CD (P_n)=
\oplus_{\alpha \in \mathbb{N}^n}P_n\der^\alpha =A_n$.
   $\Box $

By definition, the $\gm$-{\em adic topology} on the algebra $\hA$
is given by the ascending chain of {\em left ideals} of the
algebra $\hA$ (neighbourhoods  of zero)
$$ \gm^{[0]}:=\hA\supset \cdots \supset  \gm^{[i]}:= \sum_{|\alpha | \geq i} \hA\der^\alpha  \supset \cdots
\supset\cap_{i\geq 0}\gm^{[i]}=0.$$ The algebra $\hA$ is a {\em
complete} (w.r.t. the $\gm$-topology) topological algebra. The
`partial derivatives' over $A$, $D_i\in \Der_{A, c}(\hA )$,
$i=1,\ldots , s$, are {\em continuous $A$-derivations} of the
algebra $\hA$ such that
$$ D_i(\der_i)=\d_{ij},\;\;  1\leq i,j \leq s.$$

\begin{lemma}\label{i15Nov05}
For each $i=1, \ldots , s$, the map $\Psi_i(\cdot ):= \sum_{k\geq
0}(-1)^k D_i^k(\cdot ) \frac{\der_i^k}{k!}:\hA \ra \hA$, is a
homomorphism of left $\hADi$-modules where $\hADi :=
\ker_{\hA}(D_i)=A[[\der_1, \ldots , \widehat{\der_i}, \ldots ,
\der_s]]$, $ \hA = \hADi [[ \der_i]]= \hADi \oplus \hA \der_i$,
 and the following statements hold:
\begin{enumerate}
\item  the map $\Psi$ is  a projection onto the algebra  $\hADi$
of $\hA$:
$$ \Psi_i : \hA =\hADi \oplus \hA \der_i \ra \hA =\hADi \oplus \hA \der_i, \;\; a+b\der_i \mapsto a, \;\; {\rm
where}\;\; a\in \hADi, \; b\in \hA .$$ In particular,  ${\rm im}
(\Psi_i )=\hADi $ and $\Psi_i (y)=y$ for all $y\in \hADi $. \item
$\Psi_i (\der_i^k)=0$, $k\geq 1$.
\end{enumerate}
\end{lemma}

{\it Proof}. The map $\Psi_i$ is obviously well-defined
 since the algebra $\hA$ is complete and
$\der_i^k\in \gm^{[k]}$, $k\geq 0$. $\Psi_i (\der_i) =
\der_i-\der_i=0$, and $\Psi_i (y)=y$ for all $y\in A^\d$.  For any
$a\in \hA$,

\begin{eqnarray*}
 D_i\Psi_i (a)&=& D_i(a- D_i(a)\frac{\der_i}{1!}+D_i^2(a)\frac{\der_i^2}{2!}-D_i^3(a)\frac{\der_i^3}{3!}+\cdots ) \\
 &=& D_i(a)-D_i(a)-D_i^2(a)\frac{\der_i}{1!}+D_i^2(a)\frac{\der_i}{1!}+
 D_i^3(a)\frac{\der_i^2}{2!}-D_i^3(a)\frac{\der_i^2}{2!}- \cdots
 \\
 &=&0.
\end{eqnarray*}
Therefore, ${\rm im}(\Psi_i)=\hADi$.
$$ \Psi_i(\der_i^m)=\sum_{k\geq 0} (-1)^k D_i^k
(\der_i^m)\frac{\der_i^k}{k!}= (\sum_{k\geq
0}(-1)^k\frac{m(m-1)\cdots
(m-k+1)}{k!})\der_i^m=(1-1)^m\der_i^m=0.$$ Since $\hA =\hADi
[[\der_i]]$, the map $\Psi_i$ is an $\hADi$-endomorphism of the
left $\hADi$-module $\hA$ and $\Psi_i(\der_i^k)=0$, $k\geq 1$, the
map $\Psi_i$ is a projection onto the subalgebra $\hADi$ of $\hA$.
$\Box$

The map 
\begin{equation}\label{big1phis}
\Psi := \Psi_1\Psi_2\cdots \Psi_s : \hA \ra \hA, \;\;
a=\sum_{\alpha \in \mathbb{N}_s} a_\alpha \der^\alpha\mapsto \Psi
(a)=a_0
\end{equation}
is a {\em projection} onto the subalgebra $A$ of $\hA$ ($\hA =
A\oplus (\oplus_{0\neq \alpha \in  \mathbb{N}_s}A\der^\alpha$)).

\begin{theorem}\label{15Dec05}
 For any $a\in \hA := {\rm End}_K(A)$,
$$ a=\sum_{\alpha \in \mathbb{N}^s}\Psi (\frac{D^\alpha}{\alpha
!} a)\der^\alpha .$$
\end{theorem}

{\it Proof}. If $a=\sum a_\alpha \der^\alpha \in \hA $, $a_\alpha
\in A$, then, by (\ref{big1phis}), $\Psi (\frac{D^\alpha}{\alpha
!} a)=a_\alpha$. $\Box $

So, the identity map ${\rm id}_{\hA }:\hA \ra \hA$ has a nice
presentation 
\begin{equation}\label{big2phis}
{\rm id}_{\hA } (\cdot ) = \sum_{\alpha \in \mathbb{N}^s}\Psi
(\frac{D^\alpha}{\alpha !} (\cdot ))\der^\alpha .
\end{equation}

\begin{theorem}\label{s15Dec05}
 For any $\s \in \Aut_K(P_m)$,
$$\s=\sum_{\alpha \in \mathbb{N}^m} \frac{(\s (x)-x)^\alpha}{\alpha !}\der^\alpha =
{\rm id}_{P_m} +\sum_{i=1}^m (\s (x_i)-x_i)\der_i +\cdots
$$
where $\frac{(\s (x)-x)^\alpha}{\alpha !}:= \prod_{i=1}^m
\frac{(\s (x_i)-x_i)^{\alpha_i}}{\alpha_i !}$.
\end{theorem}

{\it Proof}. Let $\s'$ be the sum. Then for any $a, b \in P_m$:
\begin{eqnarray*}
\s'(ab)&=&\sum_{\alpha \in \mathbb{N}^m} \frac{(\s
(x)-x)^\alpha}{\alpha !}\der^\alpha (ab)=\sum_{\alpha \in
\mathbb{N}^m} \frac{(\s (x)-x)^\alpha}{\alpha !}\sum_{\beta +\g
=\alpha} {\alpha \choose \beta} \der^\beta  (a) \der^\g (b)\\
&=& (\sum_{\beta \in \mathbb{N}^m} \frac{(\s (x)-x)^\beta}{\beta
!}\der^\beta (a)) \, (\sum_{\g \in \mathbb{N}^m} \frac{(\s
(x)-x)^\g}{\beta !}\der^\g (b))=\s' (a) \s'(b),
\end{eqnarray*}
and so $\s'\in \Aut_K(P_m)$.  For each $i=1, \ldots , s$, $\s'
(x_i)=x_i+\s (x_i)-x_i= \s (x_i)$, hence $\s'=\s$.  $\Box$

{\it Example}. Let $\s \in \Aut_K(P_n)$, $P_n:= K[x_1, \ldots ,
x_n]$, $\s (x_i)=x_i+\l_i$ where $\l := (\l_1, \ldots , \l_n)\in
K^n$. By Theorem \ref{s15Dec05}, $$\s =\sum_{\alpha \in \Ns}
\frac{\l^\alpha \der^\alpha}{\alpha !}=\prod_{i=1}^n (\sum_{k\geq
0} \frac{(\l_i \der_i)^k}{k !})=\prod_{i=1}^n
e^{\l_i\der_i}=e^{\sum_{i=1}^n\l_i\der_i}.$$

{\it Example}. Let $\s_\l \in \Aut_K(P_n)$, $P_n:= K[x_1, \ldots ,
x_n]$, $\s (x_i)=\l_i x_i$, $\l := (\l_1, \ldots , \l_n)\in
K^{*n}$. By Theorem \ref{s15Dec05}, $\s_\l  =\sum_{\alpha \in \Ns}
(\l -1)^\alpha  \frac{x^\alpha \der^\alpha}{\alpha !}$ where $(\l
-1)^\alpha:= \prod_{i+1}^n(\l_i -1)^{\alpha_i}$. Clearly, $\s_\l
\s_\mu =\s_{\l \mu }$ for all $\l , \mu \in K^{*n}$.

Department of Pure Mathematics

University of Sheffield

Hicks Building

Sheffield S3 7RH

UK

email: v.bavula@sheffield.ac.uk

\end{document}